%% file: envelope.tex
\documentclass[letterpaper,10pt]{article}
\usepackage[textwidth=11.5cm, textheight=22cm, top=2.0cm, centering]{geometry}
\usepackage{setspace}
\usepackage[T1]{fontenc}
\usepackage{lmodern}
\usepackage[english]{babel}
\usepackage{microtype} 
\usepackage[scaled=0.95]{inconsolata} 
\makeatletter
\renewcommand{\texttt}[1]{{\ttfamily #1}}

\renewcommand{\mathtt}[1]{\text{\texttt{#1}}} 
\makeatother
\usepackage{csquotes}
\usepackage{etoolbox}
\usepackage{bbding}
\usepackage{changepage}
\usepackage{booktabs}
\usepackage{graphicx}
\usepackage[x11names, dvipsnames]{xcolor}
\definecolor{Linkz}{RGB}{30, 110, 170}
\definecolor{Darkenta}{RGB}{185, 35, 90}
\definecolor{Magentz}{RGB}{255, 35, 170}
\definecolor{Lightenta}{RGB}{254, 232, 255}
\definecolor{Reference}{RGB}{35, 180, 90}
\definecolor{Periwinkle}{RGB}{102, 51, 255}
\definecolor{yello}{RGB}{255, 245, 230}
\definecolor{Greeno}{RGB}{0, 140, 100}
\definecolor{Leeno}{RGB}{239, 255, 232}
\definecolor{Nicegreen}{RGB}{100, 200, 130}
\usepackage{amsmath, amsthm, amscd, amssymb, stmaryrd}
\usepackage{stackengine}
\usepackage{bussproofs}
\usepackage[square]{natbib}

\usepackage[
colorlinks=true,
linkcolor=Darkenta,
citecolor=Greeno,
urlcolor=Darkenta,
]{hyperref}
\usepackage{titlesec} 
\usepackage{abstract} 
\usepackage{enumitem} 
\usepackage{ccicons} 
\usepackage{float}
\usepackage{placeins} 
\usepackage{tikz} 
\usepackage{tikz-cd}
\usetikzlibrary{positioning,fit,patterns}
\usetikzlibrary{arrows.meta} 
\usepackage{adjustbox} 
\usepackage{pgfmath} 
\usepackage{ifthen} 

\newtheoremstyle{upright}
{6pt plus 2pt minus 2pt} 
{6pt plus 2pt minus 2pt} 
{\normalfont} 
{} 
{\bfseries} 
{.} 
{.5em} 
{} 
\theoremstyle{upright}
\theoremstyle{upright}
\newtheorem{theorem}{Theorem}[subsection]
\newtheorem{remark}[theorem]{Remark}

\newtheorem{thesis}[theorem]{Thesis}
\newtheorem{definition}[theorem]{Definition}

\newtheorem{proposition}[theorem]{Proposition}
\newtheorem{lemma}[theorem]{Lemma}

\newtheorem{corollary}[theorem]{Corollary}
\newtheorem{example}[theorem]{Example}

\makeatletter
\renewenvironment{proof}[1][Proof]{%
	\par\pushQED{\qed}%
	\normalfont
	\topsep6\p@\@plus6\p@\relax
	\trivlist
	\item[\hskip\labelsep\slshape #1\@addpunct{.}]%
}{%
	\popQED\endtrivlist\@endpefalse
}

\makeatother

\usepackage{algpseudocode}
\usepackage[most]{tcolorbox}
\usepackage{listings}
\newtcolorbox{breakbox}[2][]{%
	breakable,
	={#2},
	fonttitle=\bfseries,
	colback=white,
	colframe=black!20,
	coltitle=black,
	colbacktitle=white,
	boxrule=0.5pt,
	arc=0pt,
	boxsep=7pt,
	left=3pt,
	right=2pt,
	top=2pt,
	bottom=4pt,
	fontupper=\small\sffamily, 
	#1
}

\AtBeginEnvironment{quotation}{\sffamily}
\renewenvironment{quotation}
{\small\vspace{0.5em}\begin{adjustwidth}{4em}{4em}%

		\setlength{\parindent}{0pt}%
		\setlength{\parskip}{\medskipamount}%
	}
	{\end{adjustwidth}\vspace{1em}}
\newcommand{\customsectionstyle}[2]{%
	\titleformat{\section}[block]
	{\normalfont\fontsize{#1}{1.2\dimexpr#1\relax}\selectfont\centering}
	{\thesection}{1em}%
	{%
		\ifthenelse{\equal{#2}{true}}{\MakeUppercase}{\relax}%
	}%
}
\newcommand{\customsectionspacing}[3]{%
	\titlespacing*{\section}{#1}{#2}{#3}%
}
\newcommand{\customsubsectionstyle}[2]{%
	\titleformat{\subsection}[block]
	{\normalfont\fontsize{#1}{1.2\dimexpr#1\relax}\selectfont\centering}
	{\thesubsection}{1em}%
	{%
		\ifthenelse{\equal{#2}{true}}{\MakeUppercase}{\relax}%
	}%
}
\newcommand{\customsubsectionspacing}[3]{%
	\titlespacing*{\subsection}{#1}{#2}{#3}%
}
\customsectionstyle{14pt}{true}
\customsectionspacing{0pt}{3.5ex plus 1ex minus 0.2ex}{2.5ex}
\customsubsectionstyle{11pt}{true}
\customsubsectionspacing{0pt}{4ex plus 1ex minus 0.2ex}{2ex}
\usepackage{caption}

\makeatletter
\newcommand{\shorttitle}[1]{\def\@shorttitle{#1}}
\newcommand{\email}[1]{\def\@email{#1}}
\newcommand{\metadata}[1]{\def\@metadata{#1}}
\renewcommand{\maketitle}{%
	\begin{center}
		\vfill
		{\fontsize{18pt}{19pt}\selectfont \@title \par}
		\vspace{1em}
		{\normalsize \@author \par}
		\vspace{0.1em}
		{\normalsize \@date \par}
	\end{center}
}
\makeatother

\begin{document}
\input{content}
\end{document}

%% file: content.tex
\newcommand{\N}{\mathbb{N}}
\newcommand{\Pset}{\mathbb{P}}
\newcommand{\Npos}{\mathbb{N}_{>0}}
\newcommand{\NgtOne}{\mathbb{N}_{>1}}

\newcommand{\Prov}{\operatorname{Prov}}
\newcommand{\T}{\operatorname{T}}
\newcommand{\Agr}{\operatorname{Agr}}
\newcommand{\DefAgr}{\operatorname{DefAgr}}
\newcommand{\StdDefAgr}{\operatorname{StdDefAgr}}
\newcommand{\Conf}{\operatorname{Conf}}

\newcommand{\Realize}{{\vartriangleright}}
\newcommand{\Formula}{\ensuremath{\mathsf{Form}}}
\newcommand{\Bott}{\ensuremath{\bot}}
\newcommand{\Imp}{\ensuremath{\to}}
\newcommand{\Neg}[1]{\ensuremath{\lnot #1}}
\newcommand{\Eqv}[3][C]{\ensuremath{#2 \simeq_{#1} #3}}
\newcommand{\EqvRel}[1][C]{\ensuremath{\simeq_{#1}}}
\newcommand{\Regulates}[2][C]{\ensuremath{#1 \Realize #2}}
\newcommand{\LEM}[1][C]{\ensuremath{\mathsf{LEM}(#1)}}
\newcommand{\MP}[1][C]{\ensuremath{\mathsf{MP}(#1)}}
\newcommand{\Cons}[1][C]{\ensuremath{\mathsf{Cons}(#1)}}
\newcommand{\Dec}[1][C]{\ensuremath{\mathsf{Dec}(#1)}}
\newcommand{\Refute}[1][C]{\ensuremath{\mathsf{Ref}(#1)}}
\newcommand{\EvC}[1][C]{\ensuremath{\mathsf{Eval}(#1)}}
\newcommand{\NegFP}[2][C]{\ensuremath{\mathsf{NegFP}_{#1}(#2)}}
\newcommand{\evapp}[2]{\ensuremath{\mathsf{eval}(#1,#2)}}
\newcommand{\Defeq}{\ensuremath{\mathrel{\triangleq}}}

\newcommand{\Unit}{\operatorname{U}}
\newcommand{\Prime}{\operatorname{P}}
\newcommand{\Comp}{\operatorname{C}}
\newcommand{\Caught}{\operatorname{R}}


\title{AN INTUITIONISTIC\\GLANCE AT PRIMES}

\author{Milan Rosko}
\date{July 2026}

\maketitle

\begin{center}\footnotesize{
		ORCID: \href{https://orcid.org/0009-0003-1363-7158}{\footnotesize\textsf{0009-0003-1363-7158}}\\
}
\end{center}

\begin{abstract} \vspace{-0.5em}\footnotesize{
This paper gives a proof-theoretic account of how positive integers must be classified as $1$, prime, or composite in intuitionistic logic. Compositehood is expressed in $\Sigma^0_0$ by exhibiting a factorization; primality is expressed in $\Pi^0_0$ by exhibiting a lack of interior factorization. Because both searches are bounded, both predicates are decidable. Organizing the checks in stages yields a recursive sieve for the primes, a characterization of modular cancellation, and finite arithmetic certificates. The final sections distinguish what \textsc{Heyting Arithmetic} ($\mathsf{HA}$) proves internally from what depends on the standard interpretation of $\mathbb{N}$.
}
\end{abstract}

\section{Introduction}
\label{sec:introduction}

\begin{figure}[ht]
	\centering
	\includegraphics[width=.95\textwidth]{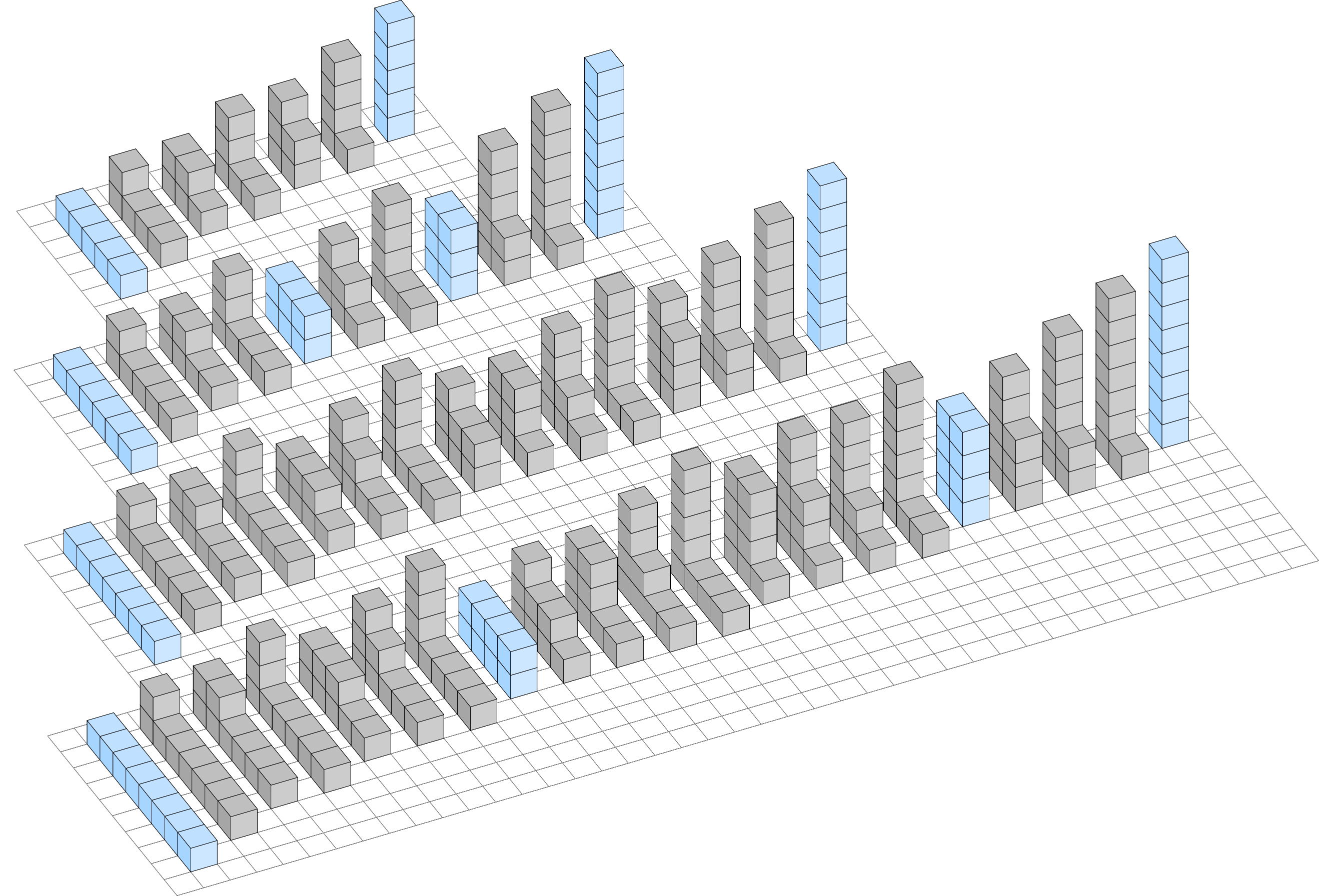}
	\caption{Within the standard bounded arithmetic presentation, $6$ and $8$ admit existential factor evidence, while organized universal refutation evidence excludes every interior factorization of $5$ and $7$. Both searches are primitive recursive; the contrast is proof polarity, not arithmetical-hierarchy complexity.}
	\label{fig:5vs6}
\end{figure}

Let mathematics correspond to rigor \emph{after} assumption; then logic is rigor \emph{before} assumption. We examine the boundary where these two intersect, from an intuitionistic perspective \citep{brouwer13,heyting31}.

Classical arithmetic proceeds as though its objects were already well formed. \textsc{Goldbach's} 1742 letter to \textsc{Euler} illustrates an older convention (of that time) that retains unity within the arithmetic of primes, with the phrase \enquote{including unity} appearing in parentheses in the original:

\begin{quotation}
	\enquote{Every integer that can be written as the sum of two primes can also be written as the sum of as many primes (including unity) as one wishes, until all terms are units.} \citet{goldbach1742letter}, in a letter to Leonhard \textsc{Euler}.
\end{quotation}

The shift in convention motivates a closer examination of the unit boundary.

Figure~\ref{fig:5vs6} makes the operative asymmetry visible before it is formalized: a successful search for $6$ or $8$ terminates with one interior factor pair, whereas the survival of $5$ or $7$ requires the exhaustion of the bounded candidate family. Accordingly, we first define three predicates through distinct evidence forms: equality $n=1$ for the unit, an interior factorization $n=ab$ with $a,b>1$ for compositehood, and a proof of $n>1$ together with refutations of every bounded interior factorization for primehood. These clauses yield the classification
\begin{equation}
	\Unit(n)\vee\Prime(n)\vee\Comp(n).
\end{equation}
Over the standard natural numbers, these clauses have their familiar extensions, and the installed bounded-search package executes the classifier for every input. The distinction is one of proof form, not numerical extension \citep{kleene45,beeson85,troelstra88}.

Second, the semantic question begins after the local classification. Its defining clauses fix the three added predicates relative to a chosen arithmetic model, while representation of the standard natural-number structure belongs to the metatheory. The paper therefore distinguishes syntactic rule conformity, internal multiplicative agreement in a fixed model, definitional agreement in a represented standard structure, whole-theory soundness, consistency, and reflection. The distinction between derivability and interpreted satisfaction is the standard syntactic--semantic distinction \citep{tarski1933}.

The arithmetic core begins with the unit boundary and the resulting ternary classifier. It then develops three consequences: staged ordering makes bounded refutation operational; finite catchers lead through \textsc{Euclidean Escape} to the sieve; and the prime--composite distinction yields modular cancellation. The dependency map in Figure~\ref{fig:arithmetic-core-dependencies} makes two structural points explicit: ternary classification is the common input to all three branches, and the catcher branch continues through escape to the sieve transition.

\begin{figure}[ht]
	\centering
	\begin{tikzpicture}[
		node distance=9mm and 10mm,
		dependency/.style={
			draw,
			rounded corners=2pt,
			align=center,
			inner xsep=15pt,
			inner ysep=5pt
		},
		arrow/.style={->, >=Stealth, semithick}
	]
		\node[dependency] (unit)
			{Unit boundary\\[-0.25em]\scriptsize\textup{Lemma~\ref{lem:unit-neither-prime-nor-composite}}};
		\node[dependency, right=of unit] (classification)
			{Ternary\\classification\\[-0.25em]
			\scriptsize\textup{Theorem~\ref{thm:ternary-classification}}};

		\node[dependency, below left=of classification] (staged)
			{Staged\\ordering\\[-0.25em]
			\scriptsize\textup{Definition~\ref{def:stage-ordering}}};
		\node[dependency, below=of classification] (catchers)
			{Finite\\catchers\\[-0.25em]
			\scriptsize\textup{Definition~\ref{def:catcher}}};
		\node[dependency, below right=of classification] (cancellation)
			{Modular\\cancellation\\[-0.25em]
			\scriptsize\textup{Theorem~\ref{thm:prime-cancellation}}};

		\node[dependency, below=of catchers] (escape)
			{\textsc{Euclidean}\\\textsc{escape}\\[-0.25em]
			\scriptsize\textup{Theorem~\ref{thm:euclidean-escape}}};
		\node[dependency, below=of escape] (sieve)
			{Sieve\\transition\\[-0.25em]
			\scriptsize\textup{Theorem~\ref{thm:least-survivor-transition}}};

		\draw[arrow] (unit) -- (classification);
		\draw[arrow] (classification) -- (staged);
		\draw[arrow] (classification) -- (catchers);
		\draw[arrow] (classification) -- (cancellation);
		\draw[arrow] (catchers) -- (escape);
		\draw[arrow] (escape) -- (sieve);
	\end{tikzpicture}
	\caption{Dependencies in the arithmetic core. The three branches from ternary classification are independent developments; the central branch continues through finite catchers, \textsc{Euclidean Escape}, and the sieve.}
	\label{fig:arithmetic-core-dependencies}
\end{figure}

The metamathematical claims connect four distinctions.
	\begin{enumerate}[label=\textnormal{(\roman*)}]
		\item Proper divisibility identifies the role of the unit in ordinary divisor semantics.
		\item Iterated defect regulation distinguishes a locally completed decision from a uniform totalization.
		\item Direct escape separates finite capture from global completeness.
		\item Definitional conservativity distinguishes interpretation within a fixed model from metatheoretic representation of $\mathbb{N}$.
	\end{enumerate}
Throughout, $\neg A$ abbreviates $A\to\bot$. Finite existential and universal searches use decidable equality, order, and multiplication on $\N$. No form of \textsc{Church's Thesis} is assumed in the arithmetic arguments \citep{church36ajm,sieg97}.

\section{Ternary Multiplicative Classification}
\label{sec:classification}

\subsection{Boundary and interior factorizations}
\label{sec:classification:support}

\begin{definition}[Boundary and interior factorizations]
	For $n\in\Npos$, let
	\begin{equation}
		\operatorname{Fac}(n) = \{(a,b)\in\Npos^2:ab=n\}
	\end{equation}
	be its set of positive multiplicative presentations. Its boundary and interior are
	\begin{equation}
	\begin{array}{c}
			\partial\operatorname{Fac}(n) = \{(a,b)\in\operatorname{Fac}(n):a=1\vee b=1\}, \\[0.5em]
			\operatorname{IntFac}(n) = \{(a,b)\in\operatorname{Fac}(n):a,b>1\}.
		\end{array}
	\end{equation}
	Hence
	\begin{equation}
		\partial\operatorname{Fac}(n)=\{(1,n),(n,1)\},
	\end{equation}
	with the two entries coinciding when $n=1$. An interior factorization is an element of $\operatorname{IntFac}(n)$.
\end{definition}

\begin{remark}
	The prime--composite distinction separates the boundary factorizations $1\cdot n$ and $n\cdot1$ from factorizations whose two factors exceed $1$. The same boundary terminates later divisibility contradictions: a prime divisor forced to divide $1$ is impossible.
\end{remark}

\begin{definition}[Arithmetic predicates]
	\label{def:arithmetic-predicates}
	The domain is $\Npos$. Write $\NgtOne$ for the positive integers strictly greater than $1$. The predicates $\Unit$, $\Prime$, and $\Comp$ are predicates on $\Npos$. The unit predicate is
	\begin{equation}
		\Unit(n)\iff n=1.
	\end{equation}
	Let the non-unit capture relation be
	\begin{equation}
		R(a,b,n) \iff a>1\wedge b>1\wedge ab=n.
	\end{equation}
	Quantification over this relation is bounded by the input. Compositehood is interior capture:
	\begin{equation}
		\Comp(n)\iff\exists a,b\le n\,R(a,b,n).
	\end{equation}
	Relative to the standard bounded-search presentation, primality is bounded refutation:
	\begin{equation}
		\Prime(n)\iff n>1\wedge \forall a,b\le n\,\neg R(a,b,n).
	\end{equation}
\end{definition}

\begin{remark}
	\label{rem:extension-regulation}
	Under the computational reading used below, unit, composite, and prime assertions use respectively equality evidence, an interior factor pair, and bounded refutation evidence.
\end{remark}

\subsection{The unit boundary}
\label{sec:classification:unit}

\begin{lemma}[The unit case]
	\label{lem:unit-neither-prime-nor-composite}
	\begin{equation}
		\Unit(1),\quad \neg\Prime(1),\quad \neg\Comp(1).
	\end{equation}
\end{lemma}

\begin{proof}
	By Definition~\ref{def:arithmetic-predicates}, $\Unit(1)$ holds because $1=1$. The predicate $\Prime(1)$ requires $1>1$, which is impossible. The predicate $\Comp(1)$ requires $R(a,b,1)$ for some $a,b$, but $R(a,b,1)$ includes $a>1$, $b>1$, and $ab=1$, which cannot hold in $\Npos$.
\end{proof}

\begin{theorem}[Ternary classification]
	\label{thm:ternary-classification}
	For every $n\in\Npos$,
	\begin{equation}
		\Unit(n)\vee\Prime(n)\vee\Comp(n).
	\end{equation}
	The three cases are mutually exclusive.
\end{theorem}

\begin{proof}
	Decide whether $n=1$. If so, $\Unit(n)$. Otherwise $n>1$. The standard bounded-pair enumeration and its total decision operations construct
	\begin{equation}
		\Comp(n)\vee\forall a,b\le n\,\neg R(a,b,n).
	\end{equation}
	The first disjunct supplies a positive interior witness. The second, together with $n>1$, is exactly $\Prime(n)$. Mutual exclusion follows because the unit case has $n=1$, whereas the other cases require $n>1$, and a witness for $R(a,b,n)$ contradicts its refutation.
\end{proof}

\begin{lemma}[Trivial divisor]
	\label{lem:no-prime-divides-unit}
	If $\Prime(q)$, then $q\nmid 1$.
\end{lemma}

\begin{proof}
	Assume $\Prime(q)$. Then $q>1$ by Definition~\ref{def:arithmetic-predicates}. If $q\mid 1$, then $1=qk$ for some $k\in\Npos$. Since $q>1$ and $k\ge1$, this gives $qk>1$, contradicting $qk=1$. Hence $q\nmid 1$.
\end{proof}

\begin{definition}[Proper divisibility in the interior semigroup]
	\label{def:strict-interior-divisibility}
	Let $D=\NgtOne$. For $a,b\in D$, define proper divisibility by
	\begin{equation}
		a\mid_D b \quad\iff\quad \exists c\in D\,(ac=b).
	\end{equation}
	Define atomicity in $D$ by
	\begin{equation}
		\operatorname{Atom}_D(p) \quad\iff\quad p\in D\wedge \neg\exists a,b\in D\,(ab=p),
	\end{equation}
	and define the prime-element law relative to proper divisibility by
	\begin{equation}
		\operatorname{PrimeEl}_D(p) \quad\iff\quad
		\begin{array}{c}
			p\in D \; \wedge\\ \forall a,b\in D\,  \bigl( p\mid_D ab \to p\mid_D a\vee p\mid_D b \bigr).
		\end{array}
	\end{equation}
	The subscript $D$ records that the cofactor, as well as the divided elements, must remain in the interior domain. Omitting the unit cofactor turns ordinary divisibility into an irreflexive proper-divisibility relation.
\end{definition}

\begin{lemma}[Proper-divisibility pathology]
\label{lem:proper-divisibility-pathology}
Let $S$ be a commutative cancellative semigroup such that $xz=x$ is impossible for all $x,z\in S$. Define
\begin{equation}
	x\prec y\quad\iff\quad\exists z\in S\,(xz=y).
\end{equation}
Then $\prec$ is irreflexive, and the prime-element condition
\begin{equation}
	\forall a,b\in S\,\bigl(p\prec ab\to p\prec a\vee p\prec b\bigr)
\end{equation}
has empty extension.
\end{lemma}

\begin{proof}
	The hypothesis gives $\neg(x\prec x)$. For $p\in S$, the factor $p$ witnesses $p\prec p^2$. The displayed implication would yield $p\prec p$, a contradiction.
\end{proof}

\begin{proposition}[Unit-removal split]
	\label{prop:unit-removal-trilemma}
	For the structure in Definition~\ref{def:strict-interior-divisibility}, proper divisibility is irreflexive:
	\begin{equation}
		\forall a\in D\,\neg(a\mid_D a);
	\end{equation}
	the atoms of the multiplicative semigroup $\NgtOne$ are exactly the ordinary prime numbers:
	\begin{equation}
		\forall p\in D\, \bigl( \operatorname{Atom}_D(p) \iff \Prime(p) \bigr);
	\end{equation}
	the prime-element predicate relative to proper divisibility is empty:
	\begin{equation}
		\forall p\in D\, \neg\operatorname{PrimeEl}_D(p).
	\end{equation}
	All three assertions are constructively provable in $\mathsf{HA}$.
\end{proposition}

\begin{proof}
	The first and third assertions are the instance $S=D$ of Lemma~\ref{lem:proper-divisibility-pathology}. For $p\in D$, $\operatorname{Atom}_D(p)$ states that every alleged factorization $p=ab$ with $a,b>1$ leads to contradiction. Such factors automatically satisfy $a,b\le p$, so this is equivalent to the bounded refutation family in Definition~\ref{def:arithmetic-predicates}. Thus $\operatorname{Atom}_D(p)\iff\Prime(p)$.
\end{proof}

\begin{remark}
	The split concerns the choice of relation, not the numerical status of ordinary primes. Atomicity retains the ordinary numerical primes. The prime-element law formed with proper divisibility has empty extension because the relation is irreflexive. Reflexively closing the relation restores ordinary divisor semantics. The reflexive closure
	\begin{equation}
		a\preceq_D b \quad\iff\quad a=b\vee a\mid_D b
	\end{equation}
	agrees with ordinary divisibility on $D$. Algebraically, the same equality case is realized by the unitization
	\begin{equation}
		D^1=D\cup\{1\}=\Npos,
	\end{equation}
	because for $a,b\in D$,
	\begin{equation}
		a\mid_{D^1}b
		\quad\iff\quad
		\exists c\in D^1\,(ac=b)
		\quad\iff\quad
		a\preceq_D b.
	\end{equation}
	The cofactor $1$ witnesses reflexivity and supplies the boundary presentations. Corollary~\ref{cor:unitization-prime-recovery} later shows that the corresponding prime-element law again selects exactly the ordinary primes.
\end{remark}

\begin{example}
	Define strict coprimality by
	\begin{equation}
		\operatorname{Coprime}_D(a,b)
		\quad\iff\quad
		\neg\exists d\in D\,
		\bigl(d\mid_D a\wedge d\mid_D b\bigr).
	\end{equation}
	Then $\operatorname{Coprime}_D(2,4)$ holds: although $2\mid_D4$, no interior element strictly divides $2$. Strict coprimality therefore classifies $2$ and $4$ differently from ordinary coprimality. Domain restriction also changes atomicity: in the multiplicative subsemigroup $\{4^k:k\ge1\}$, the ordinary composite $4$ is an atom because none of its factorizations has both factors in the selected subsemigroup.
\end{example}

\begin{remark}
	The bare semigroup
	\begin{equation}
		(\NgtOne,\;\cdot\;)
	\end{equation}
	determines its atom class. The chosen divisibility relation determines whether a prime-element law agrees with that class. Additional arithmetic vocabulary identifies the abstract atoms with the named numerals $\{2,3,5,\ldots\}$.
\end{remark}

\begin{proposition}[Failure of the binary opposition on $\Npos$]
\label{sec:classification:binary-failure}
	On $\Npos$, the implication
	\begin{equation}
		\neg\Prime(n)\to\Comp(n)
	\end{equation}
	does not hold for every $n\in\Npos$.
\end{proposition}

\begin{proof}
	Take $n=1$. By Lemma~\ref{lem:unit-neither-prime-nor-composite}, $\neg\Prime(1)$ and $\neg\Comp(1)$. Thus the extensions of $\neg\Prime$ and $\Comp$ differ at the unit.
\end{proof}

\subsection{Witness and refutation}
\label{sec:classification:witness-refutation}

\begin{lemma}[Composite witnesses]
	\label{lem:composite-witnesses}
	If $\Comp(n)$, then there exist $a,b\in\NgtOne$ such that
	\begin{equation}
		ab=n.
	\end{equation}
	Conversely, any such pair witnesses $\Comp(n)$.
\end{lemma}

\begin{proof}
	If $\Comp(n)$, then Definition~\ref{def:arithmetic-predicates} gives $a,b$ such that $R(a,b,n)$. By definition of $R$, this includes $a>1$, $b>1$, and $ab=n$, so $a,b\in\NgtOne$. Conversely, if $a,b\in\NgtOne$ and $ab=n$, then $a,b\le n$ in $\Npos$, so $R(a,b,n)$ holds and hence $\Comp(n)$ holds.
\end{proof}

\begin{lemma}[Primehood under the bounded-search interpretation]
	\label{lem:primehood-bounded-refutation}
	Under the computational interpretation of bounded quantification fixed by $\mathsf{BndPkg}(p)$ in Definition~\ref{def:bounded-search-presentation}, a proof of $\Prime(p)$ yields evidence for $p>1$ and a finite procedure producing refutations
	\begin{equation}
		R(a,b,p)\to\bot
	\end{equation}
	for every $a,b\le p$.
\end{lemma}

\begin{proof}
	By Definition~\ref{def:arithmetic-predicates}, $\Prime(p)$ includes $p>1$ and $\forall a,b\le p\,\neg R(a,b,p)$. Intuitionistically, $\neg R(a,b,p)$ is the operation $R(a,b,p)\to\bot$. The presentation $\mathsf{BndPkg}(p)$ fixes the finite index domain, enumerates its pairs, and organizes the procedure producing the dependent refutations. This bounded-universal evidence specifies the computational interpretation; individual first-order $\mathsf{HA}$ derivations may have different literal syntax.
\end{proof}

\begin{remark}
	In \textsc{Kleene's} recursive-predicate framework, the \textsc{Arithmetic Hierarchy} concerns unbounded quantifier structure; bounded quantification over recursive relations remains in the recursive base \citep{kleene52}. The characteristic functions of $\Comp$ and $\Prime$ are primitive recursive in the standard presentation. The existential/universal distinction records proof polarity, not an unbounded $\Sigma^0_1/\Pi^0_1$ or recursively enumerable/co-recursively enumerable separation.
\end{remark}

\subsection{Bounded decidability}
\label{sec:classification:decision}

\begin{proposition}[Decidability relative to a bounded presentation]
\label{prop:bounded-decidability-capture}
	For each $n\in\Npos$, the standard presentation $\mathsf{BndPkg}(n)$ makes the relation
	\begin{equation}
		R(a,b,n)
	\end{equation}
	decidable over the completed finite domain $a,b\le n$. Hence the installed presentation uniformly decides both $\Comp(n)$ and $\Prime(n)$.
\end{proposition}

\begin{proof}
	Fix $n\in\Npos$. The pair enumerator in $\mathsf{BndPkg}(n)$ lists every $(a,b)$ with $a,b\le n$. Total equality, order, and multiplication decide each instance $R(a,b,n)$. The organized finite fold therefore returns either a factor witness for $\Comp(n)$ or the complete refutation family for $\Prime(n)$.
\end{proof}

\section{Staged Block Presentations}
\label{sec:staged-search}

\subsection{Operational content of bounded notation}
\label{sec:staged-search:bounded-notation}

\begin{definition}[Bounded-search presentation]
\label{def:bounded-search-presentation}
	For a fixed input $n$, a bounded-search presentation consists of the following uniform data:
	\begin{equation}
		\mathsf{BndPkg}(n)=
		\bigl(
		D_n,e_n,\mu,\operatorname{EqDec},\operatorname{PairEnum}_n,\rho_n
		\bigr),
	\end{equation}
	where $D_n=\{0,\ldots,n\}$ is supplied as a completed finite type, $e_n$ enumerates it, $\mu:\N\times\N\to\N$ is total multiplication, $\operatorname{EqDec}$ decides equality on $\N$, $\operatorname{PairEnum}_n$ exhaustively enumerates $D_n\times D_n$, and $\rho_n$ gives a deterministic organization of that pair enumeration and its finite universal fold.
\end{definition}

\begin{proposition}[Bounded-search unpacking]
\label{prop:bounded-search-unpacking}
	An inhabitant of $\mathsf{BndPkg}(n)$ constructs the decision
	\begin{equation}
		\exists a,b\le n\,R(a,b,n)\quad+\quad\forall a,b\le n\,\neg R(a,b,n).
	\end{equation}
	The right summand contains the organized finite family of refutations determined by $\operatorname{PairEnum}_n$ and $\rho_n$.
\end{proposition}

\begin{proof}
	Enumerate every bounded pair, compute $ab$ using $\mu$, and decide $ab=n$ with $\operatorname{EqDec}$ together with the order tests in $R$. A successful test supplies the left summand. If all tests fail, the finite universal fold organized by $\rho_n$ supplies the right summand.
\end{proof}

\begin{remark}
	Although divisibility is definable by
	\begin{equation}
		d\mid n\quad\iff\quad\exists q\,(dq=n),
	\end{equation}
	its constructive decision requires a completed bounded cofactor domain, decidable arithmetic tests, exhaustive enumeration, and termination with either a quotient witness or refutations of all candidates. $\mathsf{BndPkg}(n)$ records this data. A staged block construction supplies a partial configuration, and a completion map connects it to the bounded pair square.
\end{remark}

\subsection{Raw blocks and admissible horizontal orders}
\label{sec:staged-search:blocks}

\begin{definition}[Finite block configuration]
\label{def:block-configuration}
	For an input $n\in\Npos$ and stage $s$, let $I_s(n)$ be a finite decidable index set and let
	\begin{equation}
		\beta_s^n:I_s(n)\longrightarrow\NgtOne\times\NgtOne.
	\end{equation}
	assign a factor block $\beta_s^n(i)=(a_i,b_i)$ to each index. Its decidable test is
	\begin{equation}
		E_s(n,i)\quad\iff\quad a_i b_i=n.
	\end{equation}
	Let $\prec_s^n$ be a decidable acyclic precedence relation on $I_s(n)$. The raw, vertically generated configuration is
	\begin{equation}
		V_s(n)=\bigl(I_s(n),\beta_s^n,\prec_s^n\bigr).
	\end{equation}
	Write $\mathcal C_s(n)=V_s(n)$ for the finite block configuration at stage $s$.
	Thus every block records $(a_i,b_i)$, the non-unit conditions through the codomain of $\beta_s^n$, and the test $a_i b_i\mathrel{?=}n$.
\end{definition}

\begin{definition}[Completion witness]
\label{def:block-completion-witness}
	A stage is complete for the bounded factor problem at $n$ when
	\begin{equation}
	\begin{array}{c}
		\operatorname{Complete}_s(n) \\[0.5em]
		\Updownarrow \\[0.5em]
		\forall a,b\, \Bigl[ 2\le a\le n\wedge2\le b\le n \to \exists i\in I_s(n)\,\beta_s^n(i)=(a,b)\Bigr].
		\end{array}
	\end{equation}
	This witness connects the evolving block family to the completed bounded square supplied by $\mathsf{BndPkg}(n)$.
\end{definition}

\begin{proposition}[Completed-stage recovery]
\label{prop:completion-passage-primality}
	For every $o_s$,
\begin{equation}
	\begin{array}{c}
		\operatorname{Complete}_s(n)\wedge
		\operatorname{Valid}_s(n,o_s)\\[0.5em]
		\downarrow \\[0.5em]
		\bigl(\Comp_s(n,o_s)\iff\Comp(n)\bigr)\wedge \bigl( \operatorname{Surv}_s(n,o_s)\iff\Prime(n)\bigr).
	\end{array}
\end{equation}

\end{proposition}

\begin{proof}
	If $a,b>1$ and $ab=n$, then $a,b\le n$, so every interior factor pair of $n$ lies in the bounded square. The completion witness therefore maps each such pair to a stage block, and validity places that block once in the horizontal configuration. A stage capture is an interior factor witness, while an ordinary interior factor witness becomes a stage capture through completion. This proves the composite equivalence. The same correspondence transports the complete family of factor refutations in both directions, which proves the prime equivalence.
\end{proof}

\begin{definition}[Ordering codes and positive defects]
\label{def:stage-ordering}
	An ordering code is a pair
	\begin{equation}
		o=(\pi,h),\quad\pi=(i_0,\ldots,i_{k-1}),\quad h=(B_0,\ldots,B_{\ell-1}),
	\end{equation}
	where $\pi$ proposes an index order and $h$ proposes its horizontal block representation. Write $\operatorname{Ord}_s(n,o)$ when every index is drawn from the stage set and every proposed block lies in $\NgtOne\times\NgtOne$:
	\begin{equation}
		\operatorname{Ord}_s(n,o)\quad\iff\quad\forall j<k\,(i_j\in I_s(n))\wedge\forall j<\ell\,(B_j\in\NgtOne\times\NgtOne).
	\end{equation}
	Write $[r]=\{0,\ldots,r-1\}$. The finite decidable ordering-defect type is
	\begin{equation}
		\begin{array}{l}
			W_s^{\mathrm{ord}}(n,o) = \bigl(\{0\}\times I_s(n)\bigr) \cup\bigl(\{1\}\times[k]\times[k]\bigr)\\[0.5em]
			\cup\bigl(\{2\}\times I_s(n)\times I_s(n)\times[k]\times[k]\bigr)\\[0.5em]
			\cup\{\langle3\rangle\} \cup\bigl(\{4\}\times[k]\bigr).
		\end{array}
	\end{equation}
	A code may still omit or duplicate indices. Define the decidable positive defect relation $\operatorname{BadOrd}_s(n,o,w)$ for $w\in W_s^{\mathrm{ord}}(n,o)$ by tagged witnesses of the following forms:
	\begin{equation}
		\begin{array}{l}
			w=\langle 0,i\rangle : i\in I_s(n)\text{ and }i\text{ does not occur in the index sequence }\pi,\\[0.5em]
		w=\langle 1,j,\ell\rangle : j<\ell<k\text{ and }i_j=i_\ell,\\[0.5em]
		w=\langle 2,i,j,r,t\rangle : i\prec_s^n j,\ i_r=i,\ i_t=j,\text{ and }t<r,\\[0.5em]
		w=\langle 3\rangle : k\ne\ell,\\[0.5em]
		w=\langle 4,j\rangle : j<k,\ j<\ell,\text{ and }B_j\ne\beta_s^n(i_j).
		\end{array}
	\end{equation}
	Define admissibility negatively by
		\begin{align}
			\operatorname{Valid}_s(n,o)\quad\iff\quad &\operatorname{Ord}_s(n,o)\wedge\forall w\in W_s^{\mathrm{ord}}(n,o)\\[0.5em]
			& \neg\operatorname{BadOrd}_s(n,o,w).
		\end{align}
	The corresponding horizontal configuration is
	\begin{equation}
		H_s(n,o)=h.
	\end{equation}
\end{definition}

\begin{proposition}[Defect-free codes are admissible enumerations]
\label{prop:defect-free-ordering}
	For finite $I_s(n)$, $\operatorname{Valid}_s(n,o)$ holds exactly when $o$ enumerates every stage index once and respects every dependency in $\prec_s^n$.
\end{proposition}

\begin{proof}
	The tags $0$ and $1$ exclude omission and duplication. Since every entry is typed by $\operatorname{Ord}_s(n,o)$, the resulting sequence enumerates $I_s(n)$ exactly once. Tag $2$ excludes every reversed precedence edge. Tags $3$ and $4$ establish that the proposed horizontal list has the correct length and is obtained pointwise from $\beta_s^n$.
\end{proof}

\begin{definition}[Stage-relative decisions]
\label{def:stage-relative-primality}
	The base positive capture relation is
	\begin{equation}
		\operatorname{Capture}_0(n,a,b)\quad\iff\quad a>1\wedge b>1\wedge ab=n.
	\end{equation}
	For $B=(a,b)$, abbreviate $\operatorname{Capture}_0(n,B)=\operatorname{Capture}_0(n,a,b)$. The first negative layer is
	\begin{equation}
		\operatorname{Ref}_0(\mathcal C_s,n,o)\quad\iff\quad\forall j<\lvert H_s(n,o)\rvert\,\neg\operatorname{Capture}_0\bigl(n,B_j\bigr).
	\end{equation}
	Define
	\begin{equation}
		\Comp_s(n,o)\quad\iff\quad\operatorname{Valid}_s(n,o)\wedge\exists j<\lvert H_s(n,o)\rvert\,\operatorname{Capture}_0\bigl(n,B_j\bigr),
	\end{equation}
	and
	\begin{equation}
		\operatorname{Surv}_s(n,o)\quad\iff\quad n>1\wedge\operatorname{Valid}_s(n,o)\wedge\operatorname{Ref}_0(\mathcal C_s,n,o).
	\end{equation}
	The parameter $o$ is part of the judgment because it supplies the deterministic scan order and its correctness proof.
\end{definition}

\begin{remark}
	At each stage, capture and refutation are decidable over the finite ordered block family. Their contrast is proof polarity relative to that configuration.
\end{remark}

\begin{proposition}[Deterministic local decision]
\label{prop:deterministic-stage-decision}
	For $n>1$ and every $o$ satisfying $\operatorname{Valid}_s(n,o)$, finite scanning of $H_s(n,o)$ constructs
	\begin{equation}
		\Comp_s(n,o)+\operatorname{Surv}_s(n,o).
	\end{equation}
\end{proposition}

\begin{proof}
	Each $\operatorname{Capture}_0(n,B_j)$ is decidable. Scan the finite horizontal sequence in the order supplied by $o$. The first successful test yields the existential witness required by $\Comp_s(n,o)$. If every test fails, the collected refutations yield $\operatorname{Ref}_0(\mathcal C_s,n,o)$ and hence $\operatorname{Surv}_s(n,o)$.
\end{proof}

\begin{proposition}[One-sided stage soundness]
\label{prop:stage-one-sided-soundness}
	For every stage and valid ordering,
	\begin{equation}
		\Comp_s(n,o)\to\Comp(n),\quad\operatorname{Valid}_s(n,o)\wedge\Prime(n)\to\operatorname{Surv}_s(n,o).
	\end{equation}
\end{proposition}

\begin{proof}
	A successful stage test supplies an interior factorization and hence a witness for $\Comp(n)$. A proof of $\Prime(n)$ refutes every bounded interior factorization, including every block present at stage $s$, and therefore supplies the refutations required by $\operatorname{Surv}_s(n,o)$.
\end{proof}

\begin{remark}
	The predicate $\operatorname{Surv}_s(n,o)$ records survival relative to $\mathcal C_s(n)$. Proposition~\ref{prop:completion-passage-primality} identifies it with ordinary primehood once a completion witness has been supplied. A partial stage may therefore classify a composite as a survivor.
\end{remark}

\begin{proposition}[Order invariance of the stage extension]
\label{prop:stage-order-invariance}
	If $o$ and $o'$ are valid stage-$s$ orderings, then
	\begin{equation}
		\Comp_s(n,o)\iff\Comp_s(n,o')\quad\text{and}\quad\operatorname{Surv}_s(n,o)\iff\operatorname{Surv}_s(n,o').
	\end{equation}
\end{proposition}

\begin{proof}
	Both ordering codes enumerate every member of $I_s(n)$ exactly once. Reindex an existential witness or the finite family of refutations along the induced permutation. The decision trace depends on $o$; the underlying block family determines the classified extension.
\end{proof}

\subsection{Extension and reordering obligations}
\label{sec:staged-search:extension}

\begin{definition}[Extension]
\label{def:single-block-extension}
	Write
	\begin{equation}
		\mathcal C_s(n)\rightsquigarrow\mathcal C_{s+1}(n)
	\end{equation}
	when there is an injection
	\begin{equation}
		\iota_s:I_s(n)\hookrightarrow I_{s+1}(n)
	\end{equation}
	and a unique fresh index $c_s\in I_{s+1}(n)\setminus\operatorname{im}(\iota_s)$ such that
	\begin{equation}
		\beta_{s+1}^n(\iota_s(i))=\beta_s^n(i)
	\end{equation}
	for every old index $i$. If $\beta_{s+1}^n(c_s)=(a_{c_s},b_{c_s})$, the added block is the data
	\begin{equation}
		B_{c_s}(n)=\bigl(a_{c_s},b_{c_s},E_{s+1}(n,c_s)\bigr).
	\end{equation}
	The enlarged precedence relation $\prec_{s+1}^n$ is part of the new configuration. It may place the fresh block between old blocks and may impose new dependencies among them.
\end{definition}

\begin{definition}[Defects]
\label{def:extension-validity}
	An extension specification may include a finite decidable set $K_s(n)$ of compatibility constraints and a decidable positive violation relation
	\begin{equation}
		\operatorname{Viol}_s
		(\mathcal C_s,o_s,\mathcal C_{s+1},o_{s+1},\kappa),
		\quad \kappa\in K_s(n).
	\end{equation}
	The finite decidable extension-defect type is
	\begin{equation}
		\begin{array}{l}
			W_s^{\mathrm{ext}} (\mathcal C_s,o_s,\mathcal C_{s+1},o_{s+1}) =  \\[0.5em]
			\bigl(\{0\}\times W_s^{\mathrm{ord}}(n,o_s)\bigr) \cup \bigl(\{1\}\times W_{s+1}^{\mathrm{ord}}(n,o_{s+1})\bigr) \\[0.5em]
			\cup \bigl(\{2\}\times I_s(n)\bigr) \cup \bigl(\{3\}\times I_s(n)\bigr) \cup \bigl(\{4\}\times K_s(n)\bigr).
		\end{array}
	\end{equation}
	Define $\operatorname{BadExt}_s(\mathcal C_s,o_s,\mathcal C_{s+1},o_{s+1},u)$ for $u\in W_s^{\mathrm{ext}}(\mathcal C_s,o_s,\mathcal C_{s+1},o_{s+1})$ by the following tagged defects:
	\begin{equation}
		\begin{aligned}
			& u=\langle0,w\rangle : \operatorname{BadOrd}_s(n,o_s,w),\\
			& u=\langle1,w\rangle : \operatorname{BadOrd}_{s+1}(n,o_{s+1},w),\\
			& u=\langle2,i\rangle : i\in I_s(n)\text{ and }\iota_s(i)\text{ is absent from the index component of }o_{s+1},\\
			& u=\langle3,i\rangle : \beta_{s+1}^n(\iota_s(i))\ne\beta_s^n(i),\\
			& u=\langle4,\kappa\rangle : \kappa\in K_s(n)\text{ and }\operatorname{Viol}_s(\mathcal C_s,o_s,\mathcal C_{s+1},o_{s+1},\kappa).
		\end{aligned}
	\end{equation}
	Extension admissibility is the negatively regulated condition
	\begin{equation}
		\begin{array}{c}
			\operatorname{ExtValid}_s (\mathcal C_s,o_s,\mathcal C_{s+1},o_{s+1})\\[0.5em] \Updownarrow \\[0.5em]
			\mathcal C_s(n)\rightsquigarrow\mathcal C_{s+1}(n) \wedge\operatorname{Ord}_s(n,o_s) \land\operatorname{Ord}_{s+1}(n,o_{s+1}) \\[0.5em]
			\land \forall u\in W_s^{\mathrm{ext}} (\mathcal C_s,o_s,\mathcal C_{s+1},o_{s+1})\, \neg\operatorname{BadExt}_s (\mathcal C_s,o_s,\mathcal C_{s+1},o_{s+1},u).
		\end{array}
	\end{equation}
\end{definition}

\begin{definition}[Reordering obligation]
\label{def:reordering-obligation}
	For a valid old ordering $o_s$, define
	\begin{equation}
	\begin{array}{l}
		\operatorname{Reord}_s(n,o_s) := \\[0.5em]
		\sum_{o_{s+1}} \bigl[ \operatorname{Valid}_{s+1}(n,o_{s+1}) \\[0.5em]
		\;\land\; \operatorname{ExtValid}_s (\mathcal C_s,o_s,\mathcal C_{s+1},o_{s+1}) \bigr].
	\end{array}
	\end{equation}
	No preservation of the positions or relative order of the old indices is required. A witness determines the insertion position of the new block, any displacement of old blocks, and a reconsideration of every dependency affected by $\prec_{s+1}^n$.
\end{definition}

\begin{remark}
	The extension in Definition~\ref{def:single-block-extension} preserves the old blocks but contains no operation
	\begin{equation}
		o_s\longmapsto o_{s+1}.
	\end{equation}
	A transport operation additionally requires a fixed tie-breaking order and an ordering algorithm. These choices extend the bare inclusion of block configurations with regulatory data.
\end{remark}

\begin{proposition}[Constructing a valid horizontal order]
\label{prop:construct-valid-order}
	Let $I_s(n)$ be finite, and let $\prec_s^n$ be decidable and acyclic. After fixing a decidable linear tie-breaker $\triangleleft$ on $I_s(n)$, bounded search constructs a valid ordering $o$.
\end{proposition}

\begin{proof}
	Among the remaining indices, select the $\triangleleft$-least index having no remaining $\prec_s^n$-predecessor. Acyclicity of the finite relation guarantees that such an index exists. Remove it and repeat. Strict decrease of the remaining finite set proves termination, and the selection condition proves that the resulting sequence respects every precedence edge. For the resulting
	\begin{equation}
	\pi=(i_0,\ldots,i_{k-1}),
	\end{equation}
	put
	\begin{equation}
	B_j=\beta_s^n(i_j),\quad h=(B_0,\ldots,B_{k-1}).
	\end{equation}
	Then $o=(\pi,h)$ has no defect of any of the five forms in Definition~\ref{def:stage-ordering}.
\end{proof}

\begin{remark}
	Proposition~\ref{prop:construct-valid-order} discharges the ordering-defect layer for one presented finite configuration. Extension validity additionally requires separate proofs that every compatibility constraint in $K_s(n)$ is satisfied. A normalization theorem for individual configurations therefore has a weaker input type than a constructor inhabiting $\operatorname{Reord}_s(n,o_s)$ for every extension.
\end{remark}

\subsection{Successive negative layers}
\label{sec:staged-search:negative-layers}

\begin{definition}[Stage-negative certificate]
\label{def:stage-negative-certificate}
	The order-independent negative certificate at stage $s$ is
	\begin{equation}
		\operatorname{Neg}_s(n)\quad\iff\quad\forall i\in I_s(n)\,\bigl(E_s(n,i)\to\bot\bigr).
	\end{equation}
	Its proof object is the finite dependent family of the displayed refutations.
	For every valid $o$,
	\begin{equation}
		\operatorname{Ref}_0(\mathcal C_s,n,o)\iff\operatorname{Neg}_s(n),\quad\operatorname{Surv}_s(n,o)\iff n>1\wedge\operatorname{Neg}_s(n).
	\end{equation}
\end{definition}

\begin{theorem}[Negative-layer decomposition]
\label{thm:negative-layer-decomposition}
	For a single-block extension with fresh index $c_s$, there is a constructive equivalence
	\begin{equation}
		\operatorname{Neg}_{s+1}(n)\iff\operatorname{Neg}_s(n)\wedge\bigl(E_{s+1}(n,c_s)\to\bot\bigr).
	\end{equation}
	Consequently, for valid $o_s$ and $o_{s+1}$,
	\begin{equation}
		\operatorname{Surv}_{s+1}(n,o_{s+1})\iff\operatorname{Surv}_s(n,o_s)\wedge\neg E_{s+1}(n,c_s).
	\end{equation}
\end{theorem}

\begin{proof}
	Restrict a stage-$(s+1)$ refutation family along $\iota_s$ to obtain the old family, and evaluate it at $c_s$ to obtain the new refutation. Conversely, combine the old family with the refutation at the unique fresh index. The second equivalence follows from Definition~\ref{def:stage-negative-certificate} and validity of both ordering witnesses.
\end{proof}

\begin{proposition}[Local finite flattening]
\label{prop:local-finite-flattening}
	For an inhabited valid ordering $o_s$, define the finite conjunction
	\begin{equation}
		\operatorname{Flat}_s(n,o_s)=\bigwedge_{j<\lvert H_s(n,o_s)\rvert}\neg\operatorname{Capture}_0(n,B_j).
	\end{equation}
	Then
	\begin{equation}
		\operatorname{Flat}_s(n,o_s)\iff\operatorname{Ref}_0(\mathcal C_s,n,o_s).
	\end{equation}
	If $o_{s+1}$ additionally witnesses $\operatorname{ExtValid}_s$, then, up to the permutation determined by the new horizontal order,
	\begin{equation}
		\operatorname{Flat}_{s+1}(n,o_{s+1})\iff\operatorname{Flat}_s(n,o_s)\wedge\neg\operatorname{Capture}_0\bigl(n,\beta_{s+1}^n(c_s)\bigr).
	\end{equation}
\end{proposition}

\begin{proof}
	The first equivalence is finite universal elimination and introduction along the inhabited horizontal list. For the second, extension validity identifies every old block inside the new list and identifies the unique fresh block. Reordering the finite conjunction then separates the old conjuncts from the fresh one.
\end{proof}

\begin{theorem}[Global data for uniform flattening]
\label{thm:no-uniform-flattening-from-local}
	The operation forming $\operatorname{Flat}_s(n,o_s)$ has one inhabited finite configuration as its input. Uniform flattening over a stage family additionally requires a completed stage domain $D$ and codes for the future configurations in order to inhabit a functional of type
	\begin{equation}
		\prod_{s\in D}
		\prod_{o_s}
		\left[
		\operatorname{Valid}_s(n,o_s)
		\to
		\sum_{\varphi_s}
		\bigl(
		\varphi_s\iff
		\operatorname{Ref}_0(\mathcal C_s,n,o_s)
		\bigr)
		\right]
	\end{equation}
	for every stage in a completed domain $D$. Once $D$, the configuration family, and the horizontal-list representation are supplied uniformly, the displayed flattening operation is recursive. Producing the horizontal list from a raw configuration additionally requires a normalization or ordering constructor.
\end{theorem}

\begin{proof}
	A local flattening traverses the supplied list $H_s(n,o_s)$. A function quantified over stages takes a stage domain and codes for the corresponding lists as additional inputs. Given these data, map each coded horizontal list to the conjunction of its entries. This construction is recursive. Constructing the horizontal lists from raw configurations is the separate normalization problem described in Proposition~\ref{prop:construct-valid-order}.
\end{proof}

\begin{proposition}[No upward persistence of negative certificates]
\label{prop:no-negative-persistence}
	There is no uniform implication
	\begin{equation}
		\operatorname{Neg}_s(n)\longrightarrow\operatorname{Neg}_{s+1}(n)
	\end{equation}
	over all single-block extensions.
\end{proposition}

\begin{proof}
	Take $n=25$ and a stage whose blocks draw their factors from $\{2,3\}$. Every old equality test fails, so $\operatorname{Neg}_s(25)$ holds. Extend the configuration by the block $(5,5)$. Its equality test succeeds, so $\operatorname{Neg}_{s+1}(25)$ is false.
\end{proof}

\begin{theorem}[Obligation-generating transition]
\label{thm:obligation-generating-transition}
	Given a decided stage $(\mathcal C_s(n),o_s)$ and a single-block extension, define the next-state obligation by
	\begin{equation}
		\begin{array}{l}
			\operatorname{Obl}_{s+1}(n,o_s) := \\[0.5em]
			\sum_{o_{s+1}}  \Big[
			\operatorname{Valid}_{s+1}(n,o_{s+1}) \land \operatorname{ExtValid}_s
			(\mathcal C_s,o_s,\mathcal C_{s+1},o_{s+1}) \\[0.5em]
			 \land \bigl( \Comp_{s+1}(n,o_{s+1})+\operatorname{Surv}_{s+1}(n,o_{s+1}) \bigr)\Big]
		\end{array}
	\end{equation}
	For $n>1$, there is a constructive map
	\begin{equation}
		\operatorname{Reord}_s(n,o_s)\longrightarrow\operatorname{Obl}_{s+1}(n,o_s).
	\end{equation}
	The transition therefore separates the compatibility-bearing reordering obligation from the subsequent deterministic decision.
\end{theorem}

\begin{proof}
	An inhabitant of $\operatorname{Reord}_s(n,o_s)$ supplies $o_{s+1}$ together with its ordering, validity, and extension-validity proofs. Proposition~\ref{prop:deterministic-stage-decision} decides the ordered stage. Pairing that decision with the supplied proofs inhabits $\operatorname{Obl}_{s+1}(n,o_s)$. Theorem~\ref{thm:negative-layer-decomposition} identifies the additional refutation required when the scan survives.
\end{proof}

\subsection{Iterated defect regulation}
\label{sec:staged-search:iterated-defects}

\begin{definition}[Iterated defect system]
\label{def:iterated-defect-system}
	For each $k\in\N$, let $X_k$ be a type of decision-organizing data, let $W_k(x)$ be a type of finite defect witnesses for $x\in X_k$, and let
	\begin{equation}
		\operatorname{Bad}_k(x,w),
		\quad w\in W_k(x),
	\end{equation}
	be a decidable positive defect relation. Define
	\begin{equation}
		\operatorname{Good}_k(x) \quad\iff\quad \forall w\in W_k(x)\,\neg\operatorname{Bad}_k(x,w).
	\end{equation}
	For $k\ge0$, let
	\begin{equation}
		\partial_k:X_{k+1}\longrightarrow X_k
	\end{equation}
	record the lower-level decision data organized by an object at level $k+1$. Define the iterated realization recursively by
	\begin{equation}
		\mathcal R_0(x_0)=\operatorname{Good}_0(x_0)
	\end{equation}
	and
	\begin{equation}
		\mathcal R_{k+1}(x_{k+1})=\operatorname{Good}_{k+1}(x_{k+1})\times\mathcal R_k(\partial_k x_{k+1}).
	\end{equation}
\end{definition}

\begin{remark}
	The staged prime search instantiates the schema as follows:
		\begin{equation}
			\begin{aligned}
			\operatorname{Good}_0 &: \operatorname{Ref}_0(\mathcal C_s,n,o),\\[0.5em]
			\operatorname{Good}_1 &: \forall w\in W_s^{\mathrm{ord}}(n,o)\,
			\neg\operatorname{BadOrd}_s(n,o,w),\\[0.5em]
			\operatorname{Good}_2 &: \forall u\in W_s^{\mathrm{ext}}
			(\mathcal C_s,o_s,\mathcal C_{s+1},o_{s+1})\,
			\neg\operatorname{BadExt}_s
			(\mathcal C_s,o_s,\mathcal C_{s+1},o_{s+1},u).
			\end{aligned}
		\end{equation}
	The typing and configuration clauses from Definitions~\ref{def:stage-ordering} and~\ref{def:extension-validity} complete $\operatorname{Good}_1$ to $\operatorname{Valid}_s$ and $\operatorname{Good}_2$ to $\operatorname{ExtValid}_s$. In particular,
		\begin{equation}
		\begin{array}{l}
			\operatorname{Surv}_s(n,o) \iff n>1\\[0.5em]
			\land \operatorname{Ord}_s(n,o) \; \land \; \forall w\in W_s^{\mathrm{ord}}(n,o)\,\neg\operatorname{BadOrd}_s(n,o,w) \\[0.5em]
			\land \forall B\in H_s(n,o)\,\neg\operatorname{Capture}_0(n,B),
		\end{array}
	\end{equation}
	while survival of an extension additionally requires $\operatorname{Good}_2$. Higher levels regulate the construction or extension mechanisms used at the preceding level.
\end{remark}

\begin{proposition}[Negative form and finite decidability]
\label{prop:negative-finite-layers}
	If $W_k(x)$ is finite and decidable and $\operatorname{Bad}_k$ is decidable, then $\operatorname{Good}_k(x)$ is decidable.
\end{proposition}

\begin{proof}
	Exhaustive search over $W_k(x)$ decides whether a defect witness exists. If none exists, the finite family of refutations supplies $\operatorname{Good}_k(x)$.
\end{proof}

\begin{proposition}[Generation of iterated obligations]
\label{prop:generation-iterated-obligations}
	A single-block extension changes the level-$0$ capture family and requires a new level-$1$ ordering object. Comparing the new ordering with the preceding stage introduces the level-$2$ extension object. More generally, transporting an iterated realization $\mathcal R_k$ across an extension requires a lift
	\begin{equation}
		L_k:
		\prod_{x_k\in X_k^s}
		\mathcal R_k^s(x_k)
		\longrightarrow
		\sum_{x_k'\in X_k^{s+1}}
		\mathcal R_k^{s+1}(x_k').
	\end{equation}
	The existence of $L_k$ is additional construction data. A lift at level $k$ may itself be regulated by a defect relation $\operatorname{Bad}_{k+1}$.
\end{proposition}

\begin{proof}
	At level $0$, Theorem~\ref{thm:negative-layer-decomposition} displays the fresh factor refutation. At level $1$, Definition~\ref{def:reordering-obligation} requires a new ordering witness. At level $2$, Definition~\ref{def:extension-validity} requires refutations of extension defects. Replacing these three concrete types by $X_k$, $W_k$, and $\partial_k$ gives the displayed dependent lift at every finite level.
\end{proof}

\subsection{Extension instability and availability}
\label{sec:staged-search:availability}

\begin{theorem}[Extension instability]
\label{thm:extension-instability}
	The obligations generated from an inhabitant of $\operatorname{Valid}_s(n,o_s)$ and $\operatorname{Ref}_0(\mathcal C_s,n,o_s)$ depend on the chosen single-block extension. In particular:
	\begin{enumerate}[label=\textnormal{(\roman*)}]
		\item $\operatorname{Neg}_{s+1}(n)$ additionally requires a refutation of the fresh block test;
		\item the next precedence relation determines which stage-$(s+1)$ orderings are valid;
		\item an inhabitant of $\operatorname{Reord}_s(n,o_s)$ requires an explicit constructor or separately supplied witness.
	\end{enumerate}
	More precisely, there are two one-block extensions $\mathcal C_{s+1}^{<}$ and $\mathcal C_{s+1}^{>}$ with the same old configuration, new index set, and block assignment, for which
	\begin{equation}
		\neg\exists o^+\,\Bigl(\operatorname{Valid}_{s+1}^{<}(n,o^+)\wedge\operatorname{Valid}_{s+1}^{>}(n,o^+)\Bigr).
	\end{equation}
	Here the superscripts indicate validity relative to the two respective precedence relations.
	A reordering operation takes the presented extension and its compatibility data in addition to the old ordering.
\end{theorem}

\begin{proof}
	The first claim is Proposition~\ref{prop:no-negative-persistence}. For the second, take a one-block configuration with index $i$ and add a fresh index $c$. Let $\mathcal C_{s+1}^{<}$ impose $c\prec i$, and let $\mathcal C_{s+1}^{>}$ impose $i\prec c$. Both extensions are acyclic and preserve the same old block data. Every valid two-block ordering for one reverses the sole precedence edge of the other, which proves the displayed incompatibility. The third claim follows because $\operatorname{Reord}_s(n,o_s)$ contains the new ordering and all extension-validity refutations as components.
\end{proof}

\begin{remark}
	An extension-dependent constructor may be obtained from explicit additional data: a code for $\mathcal C_{s+1}$, decidable acyclicity, a tie-breaking order, and a proof that every constraint in $K_s(n)$ has no violation. Proposition~\ref{prop:construct-valid-order} then constructs the ordering component, and the compatibility proofs complete an inhabitant of $\operatorname{Reord}_s(n,o_s)$. Extension instability concerns derivation from the old inhabitant alone.
\end{remark}

\begin{definition}[Four levels of availability]
\label{def:availability-levels}
	Fix an input $n$. The following assertions have different data:
	\begin{enumerate}[label=\textnormal{(\roman*)}]
		\item \emph{Concrete local inhabitation} at a constructed stage $s$ is
		\begin{equation}
			\sum_{o_s}\bigl(\operatorname{Valid}_s(n,o_s)\wedge\operatorname{Ref}_0(\mathcal C_s,n,o_s)\bigr).
		\end{equation}
		\item An \emph{open-ended extension rule} forms $\mathcal C_{s+1}$ from a present stage and returns the new obligation type $\operatorname{Obl}_{s+1}(n,o_s)$ without asserting that this type is inhabited.
		\item \emph{Internal pointwise existence} over a represented stage domain $D$ is
		\begin{equation}
			\forall s\in D\,\exists o_s\,\operatorname{Valid}_s(n,o_s).
		\end{equation}
		\item A \emph{uniform total selection} is a coded total function $F$ satisfying
		\begin{equation}
			\operatorname{Tot}_D(F)\wedge\forall s\in D\,\operatorname{Valid}_s(n,F(s)).
		\end{equation}
	\end{enumerate}
\end{definition}

\begin{remark}
	A specified stage-extension operation turns local inhabitation into an open-ended rule. Internal pointwise existence additionally requires a completed stage domain, a uniform coding of configurations, and proofs discharging every generated obligation. The internal quantifier ranges over the represented domain $D$. Passing from the internal formula $\forall s\exists o_s$ to a coded selector $F$ requires an appropriate choice, witnessing, or proof-extraction theorem. Under propositions-as-types with effective normalization, an inhabitant of
	\begin{equation}
		\prod_{s\in D}\sum_{o_s}\operatorname{Valid}_s(n,o_s)
	\end{equation}
	already computes such a selector; the distinction then lies between the bare first-order formula and its effective realizer.
\end{remark}

\begin{remark}
	The claims below are relative to the displayed input data. A presently inhabited finite prefix supplies the configurations and decisions already constructed. Uniform totalization additionally takes a coding of admissible continuations together with normalization, extension, and selection operations. A recursive continuation regime and compatible implementations of these operations may supply the required global data.
\end{remark}

\begin{proposition}[Prefix-only extension obstruction]
\label{thm:finite-prefix-nondetermination}
	Let
	\begin{equation}
		\operatorname{Local}_m(n)
		=
		\prod_{s<m}
		\sum_{o_s}
		\bigl(
		\operatorname{Valid}_s(n,o_s)
		\wedge
		\operatorname{Ref}_0(\mathcal C_s,n,o_s)
		\bigr)
	\end{equation}
	be the package of valid decisions for a presently constructed finite prefix. If the permitted continuations include the two opposite one-block precedence extensions of Theorem~\ref{thm:extension-instability}, then no single next-stage ordering code is valid for both continuations.
\end{proposition}

\begin{proof}
	The type $\operatorname{Local}_m(n)$ specifies configurations, orderings, and refutations for indices $s<m$. Take the common one-block prefix used in Theorem~\ref{thm:extension-instability}. One continuation imposes $c\prec i$ and the other imposes $i\prec c$. Any two-block ordering reverses one of these edges, and hence is invalid for one of the continuations. Thus the next ordering depends on the selected continuation.
\end{proof}

\begin{remark}
	A finite prefix presents the stages already constructed. A global presentation additionally comprises a coded domain of future configurations, a uniform block enumerator, a normalization algorithm, an extension operator, and a selector for further obligations. A recursive continuation regime, a tie-breaker, acyclicity proofs, and compatibility data can supply these components. In particular, a recursively coded stage family with a uniform block enumerator and compatibility solver supports a recursive extension functional. The complete bounded square in Theorem~\ref{thm:finite-stage-totalization} supplies a uniform totalization for ordinary primality.
\end{remark}

\subsection{Pre-inhabitation boundary for hardness}
\label{sec:staged-search:hardness-boundary}

\begin{definition}[Reduction datum]
\label{def:reduction-datum}
	Let $Y$ be a completed coded source instance type with predicate $A:Y\to\operatorname{Prop}$. A many-one reduction target consists of a completed coded instance type $X$ and a total predicate family
	\begin{equation}
		P:X\longrightarrow\operatorname{Prop}.
	\end{equation}
	A reduction is a total computable map $f:Y\to X$ together with
	\begin{equation}
		\forall y\in Y\,\bigl(A(y)\iff P(f(y))\bigr).
	\end{equation}
	A \textsc{Turing Reduction} analogously requires total source and target codings and a specified oracle interface for $P$. Totality of $P$ here means that $P(x)$ is a well-formed proposition for every $x\in X$; decidability is a separate property.
\end{definition}

\begin{theorem}[Pre-inhabitation boundary for hardness]
\label{thm:pre-inhabitation-hardness}
	Hardness and completeness classifications take the target instance type $X$, target predicate $P$, their codings, and the reduction operation as prior data. Typing a total reduction into an evolving block regime therefore requires an independently constructed global target presentation.
\end{theorem}

\begin{proof}
	The reduction map has codomain $X$, and its correctness field contains the expression $P(f(y))$. Both $X$ and $P$ are therefore parameters required to form the type of reduction data in Definition~\ref{def:reduction-datum}. Target presentation supplies these parameters before a term inhabiting the reduction type establishes hardness. Thus reduction theory is subsequent to target presentation.
\end{proof}

\begin{remark}
	The theorem permits hardness analysis of each fixed stage-survival predicate $\operatorname{Surv}_s(-,o_s)$ and of ordinary primality after $\mathsf{BndPkg}$ has been installed. A global target predicate ranging over future extensions additionally requires an explicit coding and continuation regime.
\end{remark}

\begin{definition}[Generative negative depth]
\label{def:generative-negative-depth}
	Write $\mathsf{NegDepth}_\omega$ for the constructive process in which every presently completed finite iterated realization
	\begin{equation}
		\mathcal R_0\triangleleft\cdots\triangleleft\mathcal R_k
	\end{equation}
	admits extensions whose realization requires a fresh negatively regulated condition $\operatorname{Good}_{k+1}$. The transition is displayed schematically as
	\begin{equation}
		\mathsf{NegDepth}_0\rightsquigarrow\mathsf{NegDepth}_1\rightsquigarrow\mathsf{NegDepth}_2\rightsquigarrow\cdots.
	\end{equation}
\end{definition}

\begin{remark}
	The notation $\mathsf{NegDepth}_\omega$ names unbounded constructive generation of negative realization layers. The index records the open-ended depth of organizational dependence, while Proposition~\ref{prop:negative-finite-layers} allows every completed finite layer to remain decidable.
\end{remark}

\begin{proposition}[Conditional properties of an evolving block regime]
\label{thm:open-ended-primality-realization}
	Let $\mathcal C_0,\mathcal C_1,\ldots$ denote successively constructed interior block configurations, with every moment represented by its presently inhabited finite prefix. Suppose each completed stage carries a locally verified deterministic ordering and its finite factor-refutation family. Then:
	\begin{enumerate}[label=\textnormal{(\roman*)}]
		\item every completed stage is locally decidable relative to its inhabited ordering data;
		\item a presented extension has the next ordering and compatibility obligation defined in Theorem~\ref{thm:obligation-generating-transition};
		\item if opposite next-stage precedence constraints are permitted, the selected continuation determines the stage-$(s+1)$ ordering;
		\item if the regime specifies a fresh defect layer at each extension, repeated extension yields the iterated conditions of Definition~\ref{def:iterated-defect-system};
		\item a global stage presentation and globally bounded search procedure require uniform stage data in addition to local finiteness;
		\item exhaustive coverage of permitted future extensions requires continuation data beyond a fixed finite configuration;
		\item hardness or completeness analysis of a global regime requires that regime's target presentation to be independently inhabited first.
	\end{enumerate}
\end{proposition}

\begin{proof}
	Clause~1 is Proposition~\ref{prop:deterministic-stage-decision}. Clause~2 is Theorem~\ref{thm:obligation-generating-transition}. Clause~3 is Proposition~\ref{thm:finite-prefix-nondetermination}. Clause~4 is Proposition~\ref{prop:generation-iterated-obligations}. Clauses~5 and~6 unpack the absence of global presentation data from a finite prefix, as explained in the following remark after Proposition~\ref{thm:finite-prefix-nondetermination}. Clause~7 is Theorem~\ref{thm:pre-inhabitation-hardness}.
\end{proof}

\begin{thesis}[Generative totality]
\label{thesis:generated-primality-realization}
	An evolving block presentation may be locally finite while lacking a supplied global continuation regime. New blocks may require a new deterministic organization of the decision objects, and a presentation that elects to regulate those organizations by further defect conditions generates the corresponding hierarchy. These are properties of the chosen presentation. Independently supplied presentation data may also support a fixed bounded formula, global ordering, or uniform selector.
\end{thesis}

\subsection{Relation to the bounded prime classifier}
\label{sec:staged-search:totalization}

\begin{theorem}[Finite totalization for a fixed input]
\label{thm:finite-stage-totalization}
	For $n\in\Npos$, let
	\begin{equation}
		I_*(n)=\{(a,b):2\le a\le n\text{ and }2\le b\le n\}.
	\end{equation}
	If $\beta_t^n$ is a bijection from $I_t(n)$ onto $I_*(n)$, then for every valid ordering $o_t$,
	\begin{equation}
		\operatorname{Surv}_t(n,o_t)\iff\Prime(n)\quad\text{and}\quad\Comp_t(n,o_t)\iff\Comp(n).
	\end{equation}
\end{theorem}

\begin{proof}
	The bijection supplies $\operatorname{Complete}_t(n)$, and validity places each indexed block exactly once in the horizontal configuration. Proposition~\ref{prop:completion-passage-primality} now gives both equivalences.
\end{proof}

\begin{thesis}[Open-ended deterministic construction]
\label{thesis:open-ended-staged-search}
	A partial block regime is deterministic once supplied with a valid horizontal order. Enlarging the raw configuration creates a new ordering obligation, and a successful reordering exposes a new decision layer. Survival at that layer requires the previous negative certificate together with a refutation of every fresh block test.

	This open-endedness concerns staged partial configurations. For each fixed $n$, the complete bounded square $I_*(n)$ totalizes the search; a fixed lexicographic order is one uniform ordering policy.
\end{thesis}

\begin{remark}
	The \textsc{Green--Tao Theorem} has the quantifier form
	\begin{equation}
		\forall k\ge1\,\exists a,d>0\,\forall i<k\,\Prime(a+id).
	\end{equation}
	An infinite prime progression would have the stronger form
	\begin{equation}
		\exists a,d>0\,\forall i\in\N\,\Prime(a+id).
	\end{equation}
	For a prime $a$, the term at $i=a$ is
	\begin{equation}
		a+ad=a(1+d),
	\end{equation}
	and is composite. Thus realization at every finite length supplies no compatible infinite prime progression. The \textsc{Green--Tao Theorem} establishes the finite realizations through a uniform ambient construction, including a transference principle and a pseudorandom majorant \citep{greentao08}. The theorem therefore illustrates how independently constructed uniform data support a global result beyond stagewise finiteness.
\end{remark}

\subsection{Complement decomposition after totalization}
\label{sec:classification:complement}

\begin{proposition}[Global complement decomposition]
\label{prop:global-complement-decomposition}
	Relative to the total bounded-search presentation of Theorem~\ref{thm:finite-stage-totalization}, every $n\in\Npos$ satisfies
	\begin{equation}
		\neg\Prime(n)\iff\Unit(n)\vee\Comp(n).
	\end{equation}
\end{proposition}

\begin{proof}
	The right-to-left implication follows from mutual exclusion. Conversely, decide $n=1$. In that case $\Unit(n)$. If $n>1$, the installed bounded enumeration and its decision operations construct $\Comp(n)\vee\forall a,b\le n\,\neg R(a,b,n)$. The second disjunct supplies $\Prime(n)$ and contradicts $\neg\Prime(n)$, so the first supplies $\Comp(n)$.
\end{proof}

\begin{proposition}[Binary recovery above the unit]
\label{prop:binary-recovery-above-unit}
	For $n\in\NgtOne$,
	\begin{equation}
		\neg\Prime(n)\iff\Comp(n).
	\end{equation}
\end{proposition}

\begin{proof}
	For $n>1$, $\Unit(n)$ is false. Proposition~\ref{prop:global-complement-decomposition} therefore reduces to $\neg\Prime(n)\iff\Comp(n)$. Its constructive content uses the totalized bounded presentation: finite search supplies either an interior witness or the complete organized refutation required for primehood.
\end{proof}

\section{Catchers and Sieves}
\label{sec:catchers}

\subsection{Finite capture stages}
\label{sec:catchers:stages}

\begin{definition}[Finite catcher]
\label{def:catcher}
A composite catcher is generated by a finite set $S\subseteq\NgtOne$. Define
\begin{equation}
	\Caught_S(n)\iff\exists d\in S\,\exists b>1\,(db=n).
\end{equation}
Equivalently, in the notation of Definition~\ref{def:strict-interior-divisibility},
\begin{equation}
	\Caught_S(n)\iff\exists d\in S\,(d\mid_D n).
\end{equation}
For $n\in\Npos$, this is also
\begin{equation}
	\Caught_S(n)\iff\exists d\in S\,(1<d<n\wedge d\mid n).
\end{equation}
It defines an infinite, decidable capture set with finitely many generators. The \textsc{Euclidean Escape Theorem} below proves that no finitely generated catcher captures every composite.
\end{definition}

\begin{proposition}[Catchers as block configurations]
\label{prop:catchers-as-block-configurations}
	For finite $S\subseteq\NgtOne$ and fixed $n$, take
	\begin{equation}
		I_S(n)=\{(d,b):d\in S\text{ and }2\le b\le n\},
	\end{equation}
	with block map $(d,b)\mapsto(d,b)$ and any valid ordering. For the induced stage configuration,
	\begin{equation}
		\Caught_S(n)\iff\Comp_s(n,o),\quad\neg\Caught_S(n)\iff\operatorname{Neg}_s(n).
	\end{equation}
	Consequently, \textsc{Euclidean Escape} constructs a composite input carrying a negative certificate for the catcher-induced stage.
\end{proposition}

\begin{proof}
	A successful block test is exactly an equation $db=n$ with $d\in S$ and $b>1$. Any such cofactor satisfies $b\le n$, so the bounded block family is extensionally complete for $\Caught_S(n)$. Exhaustiveness of a valid ordering gives both equivalences. The final statement follows from Theorem~\ref{thm:euclidean-escape}.
\end{proof}

\subsection{Soundness of finite catchers}
\label{sec:catchers:soundness}

\begin{proposition}[Monotonicity and soundness]
\label{prop:catcher-soundness}
If $S_0\subseteq S_1$, then
\begin{equation}
	\Caught_{S_0}(n)\to\Caught_{S_1}(n).
\end{equation}
For every finite $S\subseteq\NgtOne$,
\begin{equation}
	\Caught_S(n)\to\Comp(n).
\end{equation}
\end{proposition}

\begin{proof}
	A witness $d\in S_0$ is also a witness in $S_1$. For soundness, $\Caught_S(n)$ supplies $d,b>1$ with $db=n$, which is precisely an interior factor witness for $\Comp(n)$.
\end{proof}

\subsection{Unit-forcing \textsc{Euclidean Escape}}
\label{sec:catchers:escape}

\begin{remark}
	The construction below is a unit-forcing variant of the familiar \textsc{Euclidean} argument for the inexhaustibility of primes; its proof is self-contained \citep{euclidheath}.
\end{remark}

\begin{lemma}[Prime beyond a finite boundary]
	\label{lem:prime-beyond-finite-boundary}
	For every finite $S\subseteq\NgtOne$, there is a prime $p$ greater than every member of $S$.
\end{lemma}

\begin{proof}
	For $S=\varnothing$, take $p=2$. Otherwise let $M=\max S$ and $N=M!+1$. Bounded search gives a least divisor $p>1$ of $N$. It is prime, since a proper non-unit divisor of $p$ would be a smaller divisor of $N$. If $p\le M$, then $p\mid M!$ and $p\mid M!+1$, hence $p\mid1$, impossible. Thus $p>M$.
\end{proof}

\begin{theorem}[\textsc{Euclidean Escape}]
	\label{thm:euclidean-escape}
	For every finite $S\subseteq\NgtOne$, there exists $m$ such that
	\begin{equation}
	\Comp(m)\wedge\neg\Caught_S(m).
	\end{equation}
\end{theorem}

\begin{proof}
	Put
	\begin{equation}
		L=\prod_{d\in S}d,\quad m=(L+1)^2.
	\end{equation}
	The pair $(L+1,L+1)$ witnesses $\Comp(m)$, since $L+1>1$. For every $d\in S$, one has $d\mid L$ and therefore $d\mid(m-1)$, because
	\begin{equation}
		m-1=L^2+2L.
	\end{equation}
	If $d\in S$ caught $m$, then $d\mid m$ as well. Subtracting the two divisible values gives $d\mid1$, impossible because $d>1$. Thus $\neg\Caught_S(m)$.
\end{proof}

\begin{remark}
	The catcher relation uses the proper divisibility of Definition~\ref{def:strict-interior-divisibility}; the escape arguments use ordinary divisibility in $\Npos$. In both constructions, the unit is the terminal value in a divisibility contradiction: an old divisor of $L$ and $(L+1)^2$, or of $M_F$ and $M_F\pm1$, would divide $1$. The same argument applies independently at every finite stage.
\end{remark}

\begin{proposition}[Paired prime escape]
\label{prop:paired-prime-escape}
Let $F\subseteq_{\mathrm{fin}}\Pset$, with the empty product interpreted as $1$, and put
\begin{equation}
	M_F=6\prod_{p\in F}p,\quad N_F^-=M_F-1,\quad N_F^+=M_F+1.
\end{equation}
There are distinct primes $q_-,q_+\notin F$ such that
\begin{equation}
	q_-\mid N_F^-,\quad q_+\mid N_F^+.
\end{equation}
\end{proposition}

\begin{proof}
	The factor $6$ gives $N_F^-\ge5$ and $N_F^+\ge7$, including when $F$ is empty. For $p\in F$, one has $p\mid M_F$. If $p$ divided either $M_F-1$ or $M_F+1$, subtraction of the two divisible values would give $p\mid1$, contrary to Lemma~\ref{lem:no-prime-divides-unit}. Bounded search supplies a least non-unit divisor $q_-$ of $N_F^-$ and a least non-unit divisor $q_+$ of $N_F^+$. Minimality makes each divisor prime, and the preceding argument gives $q_-,q_+\notin F$. Finally, every common divisor of $N_F^-$ and $N_F^+$ divides their difference $2$. Both escape values are odd because $M_F$ is even, so their greatest common divisor is $1$. Hence $q_-\ne q_+$.
\end{proof}

\begin{definition}[Prime closure profiles]
\label{def:prime-closure-profiles}
Fix a primitive-recursive coding of finite sequences. Write $\operatorname{PrimeSeq}(s)$ when every entry of the coded sequence $s$ satisfies $\Prime$, and write $q\notin s$ when no entry equals $q$. The finite-escape sentence is
\begin{equation}
	\mathsf E_0
	\;\equiv\;
	\forall s\,
	\bigl(
	\operatorname{PrimeSeq}(s)
	\to
	\exists q\,(\Prime(q)\wedge q\notin s)
	\bigr).
\end{equation}
Two stronger global profiles are twin-prime recurrence
\begin{equation}
	\mathsf E_{\mathrm{tw}}
	\;\equiv\;
	\forall B\,\exists p\,
	\bigl(
	p>B\wedge\Prime(p)\wedge\Prime(p+2)
	\bigr),
\end{equation}
and binary \textsc{Goldbach} coverage
\begin{equation}
	\mathsf E_{\mathrm G}
	\;\equiv\;
	\forall n\,
	\bigl(
	n\ge2
	\to
	\exists p,q\,
	(\Prime(p)\wedge\Prime(q)\wedge p+q=2n)
	\bigr).
\end{equation}
In centered coordinates, the first uses $c=p+1$ and primes $c-1,c+1$. After ordering a \textsc{Goldbach} pair, elementary parity gives the symmetric form $p=n-r$, $q=n+r$ for some $r\in\N$ with $r<n$; allowing $r=0$ includes $4=2+2$.
\end{definition}

\begin{proposition}[Relations among the closure profiles]
\label{prop:closure-profile-relations}
Over elementary arithmetic,
\begin{equation}
	\mathsf E_{\mathrm{tw}}\to\mathsf E_0,\quad\mathsf E_{\mathrm G}\to\mathsf E_0.
\end{equation}
\end{proposition}

\begin{proof}
	Given a coded finite prime sequence, compute a bound $B$ for its entries. Under $\mathsf E_{\mathrm{tw}}$, choose a twin prime $p>B$; then $p$ is absent from the sequence. Under $\mathsf E_{\mathrm G}$, choose $n>\max\{B,1\}$. If both primes in a representation $p+q=2n$ were at most $B$, their sum would be at most $2B<2n$. Decidable order therefore identifies one of $p,q$ above $B$, again giving a prime absent from the sequence.
\end{proof}

\begin{remark}
	Proposition~\ref{prop:paired-prime-escape} constructively proves $\mathsf E_0$ and supplies two fresh prime divisors, although $N_F^-$ and $N_F^+$ may be composite. Reapplication yields arbitrarily long finite extensions. Twin-prime recurrence and \textsc{Goldbach} coverage require respectively uniform prime-pair construction beyond every bound and for every even total; both entail $\mathsf E_0$, which is the closure strength established by \textsc{Euclidean Escape}. After finite sequences are coded, all three are first-order arithmetical principles; an explicit reduction of the kind considered in Section~\ref{sec:semantic-agreement:iterated-closure} gives them consistency-bearing force.
\end{remark}

\begin{lemma}[Coded list form]
	\label{mech:euclid-minus-one}
	The finite-escape consequence of the minus branch of Proposition~\ref{prop:paired-prime-escape} has the list form
	\begin{equation}
		\forall l\,\bigl(\operatorname{Forall}(\Prime,l)\to\exists q\,(\Prime(q)\wedge q\notin l)\bigr).
	\end{equation}
	For the finite range $F$ of $l$, natural-number subtraction can be eliminated by recording $M_F=N_F^-+1$. Product, sequence membership, bounded least-divisor search, and the divisibility argument are primitive recursive. Eliminating $\Prime$ through Definition~\ref{def:arithmetic-predicates} therefore yields an arithmetical translation $\mathsf E_0^*$ for which $\mathsf{HA}$ proves
	\begin{equation}
	\mathsf{HA}\vdash\mathsf E_0^*.
	\end{equation}
\end{lemma}

\subsection{Limits of finite closure}
\label{sec:catchers:limits}

\begin{corollary}[No finite catcher is globally complete]
\label{cor:no-finite-catcher-complete}
For every finite $S\subseteq\NgtOne$,
\begin{equation}
	\forall n\,(\Caught_S(n)\to\Comp(n))
\end{equation}
and
\begin{equation}
	\neg\forall n\,(\Comp(n)\to\Caught_S(n)).
\end{equation}
\end{corollary}

\begin{proof}
	The first statement is Proposition~\ref{prop:catcher-soundness}. Theorem~\ref{thm:euclidean-escape} effectively constructs from $S$ a composite $m$ satisfying $\neg\Caught_S(m)$, which refutes the second universal implication.
\end{proof}

\begin{figure}[ht]
	\centering
	\includegraphics[width=0.75\textwidth]{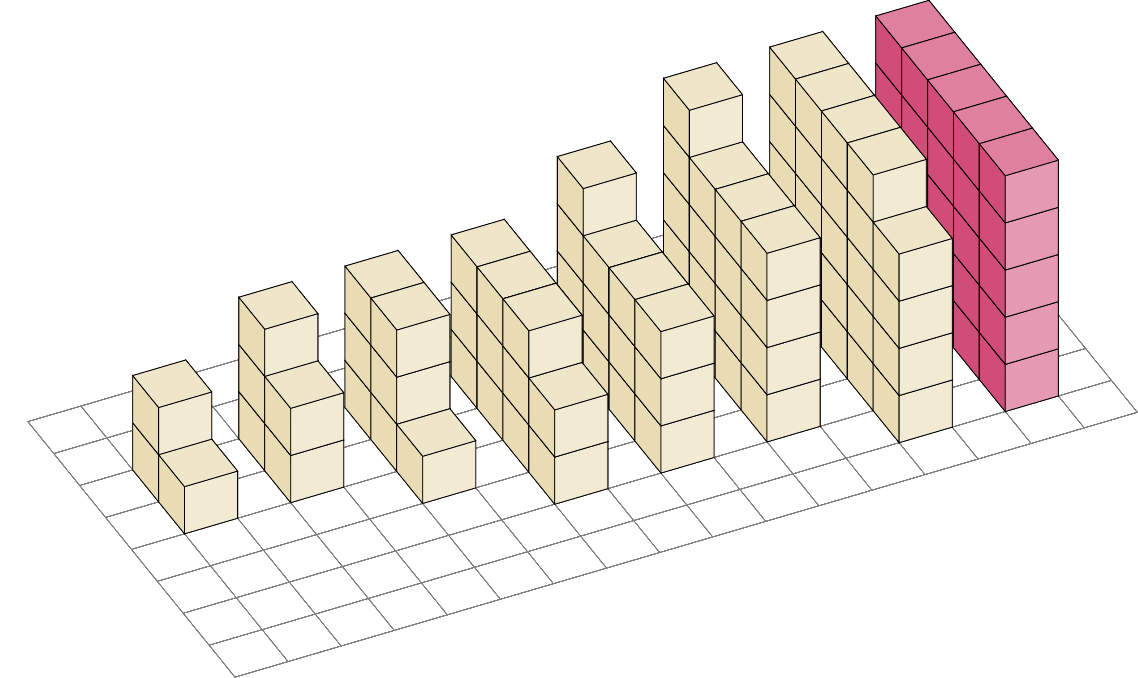}
	\caption{A catcher generated by $\{2,3\}$ would miss $25$ as composite because its factor witness $(5,5)$ lies outside the generating set.}
	\label{fig:25}
\end{figure}

Figure~\ref{fig:25} instantiates this failure for the catcher generated by $\{2,3\}$. Its survival of $25$ isolates the precise missing rule, multiplication by $5$, and thereby motivates the sieve transition in which the least survivor is promoted to a new generator.

\subsection{Sieve states}
\label{sec:sieve:states}

\begin{definition}[Recursive sieve state]
\label{def:sieve-state}
Put $P_0=\varnothing$. For a finite set $P\subseteq\NgtOne$, let
\begin{equation}
	C(P)=\bigcup_{p\in P}p\NgtOne,\quad S(P)=\NgtOne\setminus(P\cup C(P)).
\end{equation}
Given $P_k$, define
\begin{equation}
	r_k=\min S(P_k), \quad P_{k+1}=P_k\cup\{r_k\}.
\end{equation}
The existence of $r_k$ is established in Theorem~\ref{thm:least-survivor-transition}. The sets $C(P_k)$ are finitely generated capture regimes, not finite sets. Membership in each component is decidable.
\end{definition}

\begin{theorem}[Least-survivor transition]
	\label{thm:least-survivor-transition}
	For every $k$, the set $S(P_k)$ is inhabited, so the recursive clause in Definition~\ref{def:sieve-state} has a least element
\begin{equation}
	r_k=\min S(P_k).
\end{equation}
	This least survivor is prime. Moreover, $P_k$ is exactly the set of the first $k$ ordinary primes, and $r_k$ is the $(k+1)$st ordinary prime.
\end{theorem}

\begin{proof}
	Proceed by induction on $k$. The claim is immediate for $P_0$. Assume every member of $P_k$ is prime. Theorem~\ref{thm:euclidean-escape} supplies a composite $m$ with $m\notin C(P_k)$. Since no prime in $P_k$ is composite, $m\notin P_k$, hence $m\in S(P_k)$. Decidable membership gives the least survivor $r_k$.

	Suppose $r_k=ab$ with $a,b>1$. Then $a<r_k$. Minimality gives $a\in P_k\cup C(P_k)$. If $a\in P_k$, then $r_k=ab\in C(P_k)$. If $a\in C(P_k)$, write $a=qc$ with $q\in P_k$ and $c>1$; then $r_k=q(cb)\in C(P_k)$. Both cases contradict $r_k\in S(P_k)$. Thus $r_k$ is prime. Every prime smaller than $r_k$ is outside $C(P_k)$ and hence, by minimality, belongs to $P_k$. Therefore adjoining $r_k$ produces exactly the first $k+1$ primes.
\end{proof}

\begin{lemma}[Small prime divisor]
\label{lem:small-prime-divisor}
Every composite $m$ has a prime divisor $q$ such that
\begin{equation}
	q^2\le m.
\end{equation}
\end{lemma}

\begin{proof}
	Choose an interior factorization $m=ab$ with $1<a\le b$. Bounded least-divisor search gives the least divisor $q>1$ of $m$. It is prime, and $q\le a$. Hence $q^2\le a^2\le ab=m$.
\end{proof}

\subsection{Forward closure}
\label{sec:sieve:closure}

\begin{proposition}[Forward closure of sieve states]
\label{prop:sieve-forward-closure}
The sieve states form an increasing sequence of finitely generated, sound, incomplete capture regimes:
\begin{equation}
	C(P_0)\subseteq C(P_1)\subseteq C(P_2)\subseteq\cdots.
\end{equation}
At each stage, the least survivor becomes a new capture rule, while no finite stage classifies every composite.
\end{proposition}

\begin{proof}
	Since $P_k\subseteq P_{k+1}$, Definition~\ref{def:sieve-state} gives $C(P_k)\subseteq C(P_{k+1})$. Soundness follows from Proposition~\ref{prop:catcher-soundness}, and incompleteness follows from Corollary~\ref{cor:no-finite-catcher-complete}. Theorem~\ref{thm:least-survivor-transition} supplies the next prime together with its proof of non-capture.
\end{proof}

\begin{definition}[Unrestricted and interior multiplication]
\label{def:interior-multiplication-map}
Unrestricted multiplication is the map
\begin{equation}
	\mu:\Npos\times\Npos\to\Npos,\quad\mu(a,b)=ab.
\end{equation}
Its restriction to non-unit inputs is
\begin{equation}
	\mu_I:\NgtOne\times\NgtOne\to\NgtOne.
\end{equation}
\end{definition}

\begin{proposition}[Image characterization of interior multiplication]
\label{prop:relative-garden-characterization}
The map $\mu$ is surjective, while
\begin{equation}
	\operatorname{im}(\mu_I)=\{n\in\NgtOne:\Comp(n)\}.
\end{equation}
Consequently,
\begin{equation}
	\NgtOne\setminus\operatorname{im}(\mu_I)=\{p\in\NgtOne:\Prime(p)\}.
\end{equation}
\end{proposition}

\begin{proof}
	For every $n\in\Npos$, the equations $n=\mu(1,n)=\mu(n,1)$ witness surjectivity. By Definition~\ref{def:arithmetic-predicates}, membership in $\operatorname{im}(\mu_I)$ is exactly the existence of a factorization $n=ab$ with $a,b>1$, hence exactly $\Comp(n)$. Binary recovery above the unit, Proposition~\ref{prop:binary-recovery-above-unit}, identifies the remaining values with the primes.
\end{proof}

\section{Algebra}
\label{sec:algebra}

\subsection{Composite collapse}
\label{sec:algebra:collapse}

\begin{proposition}[Composite collapse]
\label{prop:composite-collapse}
If $n=ab$ with $1<a,b<n$, then in $\mathbb Z/n\mathbb Z$,
\begin{equation}
	a\not\equiv0,
	\quad
	b\not\equiv0,
	\quad
	ab\equiv0.
\end{equation}
Thus a composite presentation supplies an explicit nonzero annihilation.
\end{proposition}

\begin{proof}
	The inequalities $1<a,b<n$ make the residue classes of $a$ and $b$ nonzero. The equation $n=ab$ gives $ab\equiv0\pmod n$, witnessing failure of cancellation.
\end{proof}

\subsection{Prime cancellation}
\label{sec:algebra:cancellation}

\begin{theorem}[Prime cancellation]
	\label{thm:prime-cancellation}
	For $n>1$ and $a\in\mathbb Z/n\mathbb Z$, let $m_a(x)=ax$. Then
	\begin{equation}
	n\text{ is prime}\quad\iff\quad\forall a\not\equiv0\pmod n,\ m_a\text{ is injective}.
	\end{equation}
\end{theorem}

\begin{proof}
	Constructively, injectivity is a cancellation operation $ax=ay\Rightarrow x=y$. For prime $n$ and $a\not\equiv0\pmod n$, the \textsc{Euclidean Algorithm} computes $\gcd(a,n)=1$ and \textsc{B\'ezout Coefficients} $u,v$ with
	\begin{equation}
	ua+vn=1.
	\end{equation}
	Modulo $n$, $ua\equiv1$. Multiplying $ax\equiv ay\pmod n$ by $u$ yields $x\equiv y\pmod n$.

	Conversely, suppose every nonzero multiplication map is injective. If $n$ were composite, write $n=ab$ with $1<a,b<n$. Proposition~\ref{prop:composite-collapse} gives
	\begin{equation}
	m_a(0)=0=m_a(b)\quad\text{with}\quad0\not\equiv b\pmod n,
	\end{equation}
	contradicting injectivity. Binary recovery above the unit then gives $\Prime(n)$.
\end{proof}

\begin{corollary}[Unitized prime-element recovery]
\label{cor:unitization-prime-recovery}
For every $p\in D$, the following are equivalent:
\begin{equation}
	\Prime(p)
	\quad\iff\quad
	\forall a,b\in D^1\,
	\bigl(
	p\mid_{D^1}ab
	\to
	p\mid_{D^1}a\vee p\mid_{D^1}b
	\bigr).
\end{equation}
\end{corollary}

\begin{proof}
	Assume $\Prime(p)$ and $p\mid_{D^1}ab$. Divisibility by $p$ is decidable. If $p\mid_{D^1}a$, the first disjunct holds. Otherwise the residue of $a$ modulo $p$ is nonzero. Theorem~\ref{thm:prime-cancellation}, applied to $ab\equiv a0\pmod p$, gives $b\equiv0\pmod p$, hence $p\mid_{D^1}b$. Conversely, suppose the displayed prime-element law holds and $p=uv$ with $u,v\in D$. Then $p\mid_{D^1}uv$, witnessed by the unit cofactor, while $u,v<p$ gives $p\nmid_{D^1}u$ and $p\nmid_{D^1}v$, a contradiction. Thus $p$ has no interior factorization, and Proposition~\ref{prop:unit-removal-trilemma} gives $\Prime(p)$.
\end{proof}

\begin{theorem}[Least missed composite]
\label{thm:least-missed-composite}
The least composite not caught by the finitely generated regime $C(P_k)$ is
\begin{equation}
	r_k^2.
\end{equation}
\end{theorem}

\begin{proof}
	The pair $(r_k,r_k)$ witnesses that $r_k^2$ is composite. Suppose $p\in P_k$ caught $r_k^2$. By Corollary~\ref{cor:unitization-prime-recovery}, the prime-element law for $p$ gives $p\mid r_k$. Theorem~\ref{thm:least-survivor-transition} gives $p<r_k$, so the resulting cofactor exceeds $1$ and makes $r_k$ composite, a contradiction. Hence $r_k^2\notin C(P_k)$. If $m<r_k^2$ is composite, Lemma~\ref{lem:small-prime-divisor} supplies a prime divisor $q$ with $q^2\le m<r_k^2$, hence $q<r_k$. Theorem~\ref{thm:least-survivor-transition} gives $q\in P_k$, so $m\in C(P_k)$.
\end{proof}

\subsection{Certificates as finite proof objects}
\label{sec:algebra:certificates}

\begin{proposition}[\textsc{Lucas--Pratt} cancellation certificate]
\label{prop:lucas-pratt-cancellation}
For $n=2$, primality is verified directly. For $n>2$, a \textsc{Lucas--Pratt} certificate supplies a residue $g$, a verified equality
\begin{equation}
	n-1=\prod_{i=1}^r q_i^{e_i},
\end{equation}
and recursive primality certificates for the distinct $q_i$. Verification checks
\begin{equation}
	g^{n-1}\equiv1\pmod n
\end{equation}
and, for every $q_i$,
\begin{equation}
	g^{(n-1)/q_i}\not\equiv1\pmod n.
\end{equation}
These data construct an inverse for every nonzero residue modulo $n$ and hence certify $\Prime(n)$.
\end{proposition}

\begin{proof}
	The equation $g^{n-1}\equiv1\pmod n$ makes $g$ a unit, so its order divides $n-1$. For every distinct prime divisor $q_i$ of $n-1$, the test $g^{(n-1)/q_i}\not\equiv1\pmod n$ forces the order to contain the full $q_i$-power occurring in $n-1$. Hence $\operatorname{ord}_n(g)=n-1$. The powers $1,g,\dots,g^{n-2}$ are consequently $n-1$ distinct nonzero residues. Since there are exactly $n-1$ nonzero residue classes, they exhaust those classes. Given a nonzero residue $a$, bounded search finds $r<n-1$ with $a=g^r$, and
	\begin{equation}
	a^{-1}=g^{n-1-r}
	\end{equation}
	constructs its inverse. Every nonzero modular action is therefore cancellative, so Theorem~\ref{thm:prime-cancellation} gives $\Prime(n)$. Each smaller prime in the factorization enters through its recursive subcertificate.
\end{proof}

\subsection{Factorization as coordinate normal form}
\label{sec:coordinates:factorization}

\begin{proposition}[Terminating prime factorization]
\label{prop:terminating-coordinate-normalization}
For every $n\in\Npos$, bounded least-divisor search constructs a finite list of primes with product $n$.
\end{proposition}

\begin{proof}
	For $n>1$, bounded search finds the least divisor $\ell(n)>1$. It is prime, since a proper non-unit divisor of $\ell(n)$ would be a smaller divisor of $n$. If $\ell(n)<n$, compute the cofactor $c(n)$ with $\ell(n)c(n)=n$ and recurse on $c(n)<n$. This strict descent terminates at prime leaves. For $n=1$, return the empty list.
\end{proof}

\begin{theorem}[Uniqueness of prime factorization]
\label{thm:unique-prime-factorization}
Every finite prime factorization of $n\in\Npos$ has the same multiplicity for each prime.
\end{theorem}

\begin{proof}
	Use induction on the length of one factorization. If $p$ is its first prime and $p\mid q_1\cdots q_r$, repeated application of Corollary~\ref{cor:unitization-prime-recovery} gives $p\mid q_i$ for some $i$. Since $q_i$ is prime, the same corollary and the absence of interior factors force $p=q_i$. Reorder the second factorization to place $q_i=p$ first, then apply the ordinary cancellation law for positive natural numbers, available in $\mathsf{HA}$. Continue inductively.
\end{proof}

\subsection{Prime-coordinate support}
\label{sec:coordinates:support}

\begin{theorem}[Prime-coordinate representation]
\label{thm:prime-coordinate-representation}
Every $n\in\Npos$ has a unique finite-support exponent vector
\begin{equation}
	(v_p(n))_{p\in\Pset}\in\N^{(\Pset)}
\end{equation}
such that
\begin{equation}
	n=\prod_{p\in\Pset}p^{v_p(n)}.
\end{equation}
Thus the multiplicative monoid has coordinate normal form
\begin{equation}
	\Npos\cong\N^{(\Pset)}.
\end{equation}
\end{theorem}

\begin{proof}
	Proposition~\ref{prop:terminating-coordinate-normalization} supplies a finite prime factorization. Theorem~\ref{thm:unique-prime-factorization} makes the exponent assigned to each prime independent of the factorization. Multiplication adds exponent vectors, so the correspondence is a monoid isomorphism.
\end{proof}

\begin{corollary}[Classification in coordinates]
\label{cor:coordinate-classification}
The unit is the zero vector, a prime $p$ is the basis vector $e_p$, and composites are precisely the finite-support vectors with
\begin{equation}
	\sum_{p\in\Pset}v_p(n)\ge2.
\end{equation}
\end{corollary}

\begin{proof}
	The empty factorization represents $1$, a single prime factor represents a basis vector, and an interior factorization exists exactly when the total multiplicity is at least two.
\end{proof}

\section{Semanticity}
\label{sec:semantic-agreement}

\subsection{The stabilized formal theory}
\label{sec:semantic-agreement:theory}

\begin{definition}[Stabilized prime theory]
\label{def:stabilized-prime-theory}
Let $\Theta_{\mathrm{U,P,C}}$ collect the regulated defining clauses
\begin{equation}
	\Unit(n)\iff n=1,
\end{equation}
where $R$ is the arithmetic abbreviation
\begin{equation}
	R(a,b,n)\iff a>1\wedge b>1\wedge ab=n,
\end{equation}
and let the regulated non-unit clauses be
\begin{equation}
	\Comp(n)\iff\exists a,b\le n\,R(a,b,n),
\end{equation}
and
\begin{equation}
	\Prime(n)\iff n>1\wedge\forall a,b\le n\,\neg R(a,b,n).
\end{equation}
The local theory is the definitional extension
\begin{equation}
	\T_{\Pset}^{\mathrm{loc}}=\mathsf{HA}+\Theta_{\mathrm{U,P,C}}.
\end{equation}
For a recursively enumerable family $\Gamma$ of further principles governing global prime behavior, put
\begin{equation}
	\T_{\Pset}(\Gamma)=\T_{\Pset}^{\mathrm{loc}}+\Gamma.
\end{equation}
When $\Gamma$ is fixed, write $\T_{\Pset}$ for this extension.
\end{definition}

\begin{theorem}[Definitional conservativity and unique expansion]
\label{thm:definitional-conservativity}
There is an effective translation $A\mapsto A^{\circ}$ from the language of $\T_{\Pset}^{\mathrm{loc}}$ to the language of $\mathsf{HA}$ which replaces $\Unit$, $\Prime$, and $\Comp$ by their defining arithmetic formulas. The local theory proves
\begin{equation}
	\T_{\Pset}^{\mathrm{loc}}\vdash A\iff A^{\circ}.
\end{equation}
Consequently, for every sentence $\varphi$ in the original language of $\mathsf{HA}$,
\begin{equation}
	\T_{\Pset}^{\mathrm{loc}}\vdash\varphi\quad\iff\quad\mathsf{HA}\vdash\varphi.
\end{equation}
Moreover, every classical first-order model $\mathcal M\models\mathsf{HA}$, under ordinary \textsc{Tarskian Satisfaction}, has a unique expansion to the extended language satisfying $\Theta_{\mathrm{U,P,C}}$.
\end{theorem}

\begin{proof}
	Define $A^{\circ}$ recursively on formulas, replacing each occurrence of an added predicate by its displayed defining formula. Induction on $A$ gives the first equivalence from the defining biconditionals. If $\varphi$ is an arithmetic sentence, then $\varphi^{\circ}=\varphi$; conservativity follows by eliminating definitions from a $\T_{\Pset}^{\mathrm{loc}}$-derivation. In a fixed classical first-order arithmetic model, each defining biconditional fixes the extension of its added predicate, which gives existence and uniqueness of the expansion.
\end{proof}

\begin{remark}
	The defining clauses constrain the extensions once a domain and its arithmetic operations are interpreted. Under the computational interpretation used in this paper, $\Unit(n)$ has equality evidence, $\Comp(n)$ has an interior factor witness, and $\Prime(n)$ has evidence for $n>1$ together with the bounded refutations specified in Definition~\ref{def:arithmetic-predicates}. Theorem~\ref{thm:definitional-conservativity} shows that this local enrichment adds no theorems in the original language of $\mathsf{HA}$.
\end{remark}

\subsection{Three levels of certification}
\label{sec:semantic-agreement:levels}

\begin{definition}[Rule conformity]
\label{def:rule-conformity}
For an effectively presented theory $\T$ with decidable proof checking, let $\sigma\in\operatorname{Sent}(\T)$ be the \textsc{G\"odel} code of a $\T$-sentence and define
\begin{equation}
	\Conf_{\T}(d,\sigma)
	\;\equiv\;
	\operatorname{Prf}_{\T}(d,\sigma).
\end{equation}
Thus
\begin{equation}
	\Conf_{\T}\bigl(d,\ulcorner\Prime(\overline 5)\urcorner\bigr)
\end{equation}
states that $d$ is accepted by the rules of $\T$ as a proof of $\Prime(5)$. A canonical prime certificate contains evidence for $5>1$ and the required bounded refutations; its checker verifies that the certificate induces such a $\T$-derivation.
For a recursively enumerable extension $\T_{\Pset}(\Gamma)$, occurrences of axioms from $\Gamma$ carry enumeration-stage witnesses, preserving decidable checking of complete proof codes.
\end{definition}

\begin{remark}
	Rule conformity is syntactic. It establishes that a finite object has the shape and dependencies required by the installed calculus. Arithmetic structures, satisfaction relations, and the intended denotation of $\Prime$ enter at the adequacy level.
\end{remark}

\begin{definition}[Internal multiplicative agreement]
\label{def:agreement-with-primes}
Work in a metatheory $S$. In this section, an arithmetic structure is a classical first-order structure for the language of $\mathsf{HA}$, with satisfaction understood in the ordinary \textsc{Tarskian} sense. Let $\mathcal M$ be such a structure represented in $S$. Let $I$ fix the interpretations of $0$, $1$, addition, multiplication, equality, and order in $\mathcal M$, and interpret the added predicates $\Unit$, $\Prime$, and $\Comp$. With the arithmetic operations evaluated in $\mathcal M$, define
\begin{equation}
	\Unit_{\mathcal M}^{\times}(x)\;\equiv\;x=1,
\end{equation}
\begin{equation}
	\Comp_{\mathcal M}^{\times}(x)\;\equiv\;\exists a,b\in|\mathcal M|\,(a>1\wedge b>1\wedge ab=x),
\end{equation}
and
\begin{equation}
	\Prime_{\mathcal M}^{\times}(x)
	\;\equiv\;
	x>1\wedge
	\forall a,b\in|\mathcal M|\,
	(a>1\wedge b>1\to ab\ne x).
\end{equation}
The internal multiplicative-agreement relation is the semantic assertion
\begin{equation}
	\begin{aligned}
	\DefAgr_{\mathcal M,I}(\Theta_{\mathrm{U,P,C}})
	\;\equiv\;
	\forall x\in|\mathcal M|\,
	\bigl[&
	\bigl(
	I(\Unit)(x)
	\iff
	\Unit_{\mathcal M}^{\times}(x)
	\bigr)\\
	{}\wedge{}&
	\bigl(
	I(\Prime)(x)
	\iff
	\Prime_{\mathcal M}^{\times}(x)
	\bigr)\\
	{}\wedge{}&
	\bigl(
	I(\Comp)(x)
	\iff
	\Comp_{\mathcal M}^{\times}(x)
	\bigr)
	\bigr].
	\end{aligned}
\end{equation}
It compares the interpreted extensions of the three regulated predicates with multiplicative relations in the already fixed structure $\mathcal M$. For a standard structure $\N$ represented in $S$, write standard definitional agreement as
\begin{equation}
	\StdDefAgr_I(\Theta_{\mathrm{U,P,C}})\quad\equiv\quad\DefAgr_{\N,I}(\Theta_{\mathrm{U,P,C}}).
\end{equation}
\end{definition}

\begin{remark}
	$\mathsf{HA}$ proves that $a,b>1\wedge ab=x$ implies $a,b\le x$, so the bounded defining clauses coincide with the displayed unbounded multiplicative clauses in every classical model of $\mathsf{HA}$. For the unique definitional interpretation $I_{\mathrm{def}}$, representation of the standard structure $\N$ in the metatheory makes $\StdDefAgr_{I_{\mathrm{def}}}(\Theta_{\mathrm{U,P,C}})$ an immediate consequence of the defining clauses. The remaining metatheoretic issue is the status of $\N$ as the intended structure: first-order satisfaction determines the three predicate extensions relative to a model, while selection of the standard model requires additional representational data.
\end{remark}

\begin{definition}[Soundness, consistency, and reflection]
\label{def:soundness-consistency-reflection}
Work either in an external set-theoretic metatheory $S$, where satisfaction for all $\T$-sentences is defined, or fix a syntactic class $\Delta$ for which $S$ supplies a partial satisfaction predicate. In the latter case, $\sigma\in\Delta$ is added to the antecedent below and $\operatorname{Sat}$ is read as $\operatorname{Sat}^{\Delta}$. For the sentences in scope, whole-theory soundness has the form
\begin{equation}
	\operatorname{Sound}_{\mathcal M,I}^{S}(\T)
	\;\equiv\;
	\forall d,\sigma\,
	\bigl(
	\bigl(
	\sigma\in\operatorname{Sent}(\T)
	\wedge
	\operatorname{Prf}_{\T}(d,\sigma)
	\bigr)
	\to
	\operatorname{Sat}_{\mathcal M,I}(\sigma)
	\bigr).
\end{equation}
The provability predicate is induced by the proof relation:
\begin{equation}
	\Prov_{\T}(\sigma)\;\equiv\;\exists d\,\operatorname{Prf}_{\T}(d,\sigma).
\end{equation}
The standard consistency assertion is
\begin{equation}
	\operatorname{Con}(\T)\;\equiv\;\forall d\,\neg\operatorname{Prf}_{\T}(d,\ulcorner\bot\urcorner).
\end{equation}
For a specified formula class $\Delta$, reflection is the schema
\begin{equation}
	\operatorname{RFN}_{\Delta}(\T)
	=
	\left\{
	\forall\vec x\,
	\bigl(
	\Prov_{\T}(\ulcorner\varphi(\dot{\vec x})\urcorner)
	\to\varphi(\vec x)
	\bigr)
	:
	\varphi\in\Delta
	\right\}.
\end{equation}
\end{definition}

\begin{proposition}[Separation of certification levels]
\label{prop:certification-levels}
Rule conformity, internal multiplicative agreement, standard definitional agreement, and reflective certification have different data and logical roles. In particular:
	\begin{enumerate}[label=\textnormal{(\roman*)}]
	\item $\Conf_{\T}(d,\sigma)$ is a mechanically checkable relation on finite syntax;
	\item $\DefAgr_{\mathcal M,I}(\Theta_{\mathrm{U,P,C}})$ records the definitional interpretation inside a fixed model, while $\StdDefAgr_I$ evaluates it in a represented $\N$;
	\item $\operatorname{Con}(\T)$ and $\operatorname{RFN}_{\Delta}(\T)$ quantify over the derivations of $\T$.
\end{enumerate}
Internal multiplicative agreement certifies the defining predicates inside a fixed model. Certification of every rule of $\T$ requires whole-theory soundness, while consistency ranges over proof codes and leaves the representation of the standard structure as separate data. Full soundness supplies a bridge:
\begin{equation}
	\operatorname{Sound}_{\mathcal M,I}^{S}(\T)
	\wedge
	\neg\operatorname{Sat}_{\mathcal M,I}
	\bigl(\ulcorner\bot\urcorner\bigr)
	\longrightarrow
	\operatorname{Con}(\T).
\end{equation}
\end{proposition}

\begin{proof}
	The first three claims unpack Definitions~\ref{def:rule-conformity},~\ref{def:agreement-with-primes}, and~\ref{def:soundness-consistency-reflection}. For the final implication, a proof of $\bot$ would, by soundness, make the interpretation of $\bot$ satisfied in $\mathcal M$, contradicting the second conjunct.
\end{proof}

\subsection{Relative representation}
\label{sec:semantic-agreement:relative-adequacy}

\begin{thesis}[Relative representation of standard agreement]
\label{thesis:relative-adequacy}
Let $S$ represent the syntax of $\T$, the standard arithmetic structure $\N$, and an interpretation $I$. In the present framework, standard definitional agreement has the relative form
\begin{equation}
	S\vdash\StdDefAgr_I(\Theta_{\mathrm{U,P,C}}).
\end{equation}
For $I_{\mathrm{def}}$, the derivation records the defining clauses within the represented $\N$; representation of $\N$ is a prior metatheoretic datum.
\end{thesis}

\begin{remark}
	\label{rem:verifier-sequence}
	Let $P_0$ be the formal prime package, and let $P_1$ verify that $P_0$ satisfies a formalized agreement contract. The certificate produced by $P_1$ establishes conformity to that contract. Treating the contract as the correct interpretation of primality invokes a metatheory and an interpretation of its vocabulary. A further verifier $P_2$ can check the derivation produced by $P_1$, but its certificate is read within a further accepted proof practice. A terminating verifier sequence is accepted relative to a metatheory, interpretation, and checking practice. Additional verification relocates the point of semantic commitment, so semantic adequacy remains relative to an interpreted proof practice.

	The verifier sequence is a relativity claim. Formal checking remains exact and useful within the framework whose rules and interpretation are accepted. Within the present distinction, additional internal rules establish further conformity, while adequacy remains a relation between the package and an interpretation. If $\T$ proves a sentence named \enquote{agreement}, its proof object establishes rule conformity inside $\T$; reading that sentence as agreement with primes still invokes an interpretation of its vocabulary.
\end{remark}

\subsection{Nonstandard interpretations}
\label{sec:semantic-agreement:nonstandard}

\begin{proposition}[Nonstandard models satisfy the formal clauses]
\label{prop:nonstandard-prime-models}
Let $\mathcal M$ be a nonstandard classical first-order model of \textsc{Peano Arithmetic}, hence a classical model of every theorem of $\mathsf{HA}$. Expand $\mathcal M$ by interpreting $\Unit$, $\Prime$, and $\Comp$ through $\Theta_{\mathrm{U,P,C}}$. The expansion satisfies the defining clauses and the internally stated ternary-classification theorem, although its arithmetic domain is not isomorphic to $\N$.
\end{proposition}

\begin{proof}
	The new predicates are defined from equality, order, and multiplication already interpreted in $\mathcal M$, so the defining biconditionals hold in the expansion. The proof of ternary classification is formalizable in the arithmetic theory and therefore holds internally. Standard compactness arguments supply nonstandard models of first-order arithmetic \citep{boolos07,kaye91,hajekpudlak93}. The definitional expansion preserves their nonstandard domains.
\end{proof}

\begin{remark}
	For a nonstandard element $c\in\mathcal M$, the initial segment
	\begin{equation}
		\{x\in|\mathcal M|:x\le c\}
	\end{equation}
	is internally bounded and is treated by $\mathcal M$ as finite, while it may be externally infinite. The internally defined prime predicate expresses multiplicative primality relative to $\mathcal M$. The lesson is non-categoricity: satisfaction of the formal clauses selects the prime relation within a model, and selection of the standard model remains a separate metatheoretic task.
\end{remark}

\subsection{Philosophical readings}
\label{sec:semantic-agreement:readings}

\begin{remark}
	The arithmetic theorems establish the regulated classifier and its familiar numerical extension over the standard natural numbers. The preceding definitions support the relative-adequacy thesis by distinguishing syntactic conformity, interpreted predicate agreement, and reflective certification. A choice among realism, structuralism, inferentialism, and antirealism requires further philosophical premises.

	A realist takes the intended arithmetic structure to answer to mind-independent numerical facts and regards adequacy as correspondence with them. A structuralist locates adequacy in preservation across the appropriate arithmetic structures and mappings. An inferentialist takes the regulated assertion and transformation rules to constitute the use, and hence the meaning, of the prime predicate. The antirealist reading pursued here identifies semantic entitlement with warranted assertion inside an accepted proof practice and treats stronger adequacy proofs as relative extensions of that practice. The arithmetic theorems and metatheoretic distinctions provide common ground; none entails the antirealist reading \citep{dummett75,dummett01}.
\end{remark}

\begin{definition}
\label{def:fixed-behavioral-contract}
After $S$, $\mathcal M$, $I$, and a coding of observable behavior have been fixed, let $\mathcal C(f)$ assert that the input--output graph of a partial computable function $f$ satisfies the selected prime-semantic contract. Assume that $\mathcal C$ is extensional,
\begin{equation}
	f=g\longrightarrow\bigl(\mathcal C(f)\iff\mathcal C(g)\bigr),
\end{equation}
and nontrivial: some partial computable functions satisfy it and some do not.
\end{definition}

\begin{proposition}[\textsc{Rice} boundary for a fixed contract]
\label{prop:rice-fixed-contract}
Let $(\varphi_e)_{e\in\N}$ effectively enumerate the partial computable functions. The index set
\begin{equation}
	I_{\mathcal C}=\{e\in\N:\mathcal C(\varphi_e)\}
\end{equation}
is undecidable. Hence there is no effective inhabitant of
\begin{equation}
	\prod_{e\in\N}
	\bigl(
	\mathcal C(\varphi_e)
	+
	\neg\mathcal C(\varphi_e)
	\bigr).
\end{equation}
\end{proposition}

\begin{proof}
	\textsc{Rice's Theorem} makes every nontrivial extensional property of partial computable functions undecidable from program indices \citep{rice53}. An effective inhabitant of the displayed type would expose one of its summands for every $e$ and therefore compute a total characteristic function for $I_{\mathcal C}$.
\end{proof}

\begin{remark}
	After the contract $\mathcal C$ has been fixed and justified, \textsc{Rice's Theorem} limits its uniform total decision for arbitrary machine indices. A separate theorem in the chosen metatheory connects machine behavior to $\StdDefAgr_I$. Proposition~\ref{prop:rice-fixed-contract} permits verification of a fixed canonical implementation, proof-carrying certificates, sound incomplete verification procedures, and relative definitional-agreement proofs in a chosen metatheory.
\end{remark}

\subsection{Reflection and the \textsc{G\"odel} boundary}
\label{sec:semantic-agreement:godel}

\begin{theorem}[Conditional internal reflection boundary]
	\label{thm:conditional-internal-reflection}
	Let $\T$ be a consistent, recursively enumerable extension of $\mathsf{HA}$, and let $\Prov_{\T}$ be its standard arithmetized provability predicate satisfying the \textsc{Hilbert--Bernays--L\"ob Derivability Conditions}. Then
	\begin{equation}
	\T\nvdash\operatorname{Con}(\T).
	\end{equation}
\end{theorem}

\begin{proof}
	This is \textsc{G\"odel's Second Incompleteness Theorem} in its standard arithmetized form \citep{godel31,loeb55,hajekpudlak93}.
\end{proof}

\begin{remark}
	Theorem~\ref{thm:conditional-internal-reflection} concerns standard syntactic consistency. The formula $\operatorname{Con}(\T)$ codes derivability, whereas $\DefAgr$ fixes definitions in a model and $\StdDefAgr$ is formulated after the standard interpretation has been represented in a metatheory. Applying the theorem to a soundness or agreement principle requires an explicit implication from that principle to $\operatorname{Con}(\T)$; under the stated hypotheses, the theorem then obstructs the stronger reflective principle. A consistent theory may admit unintended models, and correct classification of every standard numeral may coexist with failure of broader reflection. A metatheory may establish definitional agreement, soundness, or consistency as relative results.
\end{remark}

\subsection{Closure principles under consistency reductions}
\label{sec:semantic-agreement:iterated-closure}

\begin{proposition}[Consistency-reduction obstruction]
\label{prop:consistency-reduction-obstruction}
Let $\T$ satisfy the hypotheses of Theorem~\ref{thm:conditional-internal-reflection}, and let $C$ be a closed sentence of the recursively presented language of $\T$. If
\begin{equation}
	\T\vdash\bigl(C\to\operatorname{Con}(\T)\bigr),
\end{equation}
then
\begin{equation}
	\T\nvdash C.
\end{equation}
\end{proposition}

\begin{proof}
	If $\T\vdash C$, modus ponens with the assumed reduction would give $\T\vdash\operatorname{Con}(\T)$, contrary to Theorem~\ref{thm:conditional-internal-reflection}.
\end{proof}

\begin{corollary}[Prime escape does not yield internal consistency]
\label{cor:escape-does-not-yield-consistency}
Assume that $\T_{\Pset}^{\mathrm{loc}}$ is consistent. Then
\begin{equation}
	\T_{\Pset}^{\mathrm{loc}}\nvdash\bigl(\mathsf E_0\to\operatorname{Con}(\T_{\Pset}^{\mathrm{loc}})\bigr).
\end{equation}
\end{corollary}

\begin{proof}
	Lemma~\ref{mech:euclid-minus-one} gives a $\mathsf{HA}$ proof of the arithmetical translation $\mathsf E_0^*$. Since $\T_{\Pset}^{\mathrm{loc}}$ is the corresponding recursively enumerable definitional extension of $\mathsf{HA}$, it proves
	\begin{equation}
	\T_{\Pset}^{\mathrm{loc}}\vdash\mathsf E_0.
	\end{equation}
	If it also proved the displayed implication, it would prove its own consistency, contradicting Theorem~\ref{thm:conditional-internal-reflection}.
\end{proof}

\begin{remark}
	Corollary~\ref{cor:escape-does-not-yield-consistency} specializes \textsc{G\"odel's Second Incompleteness Theorem}. \textsc{Euclidean} arithmetic supplies the already-provable antecedent $\mathsf E_0$; the implication $\mathsf E_0\to\operatorname{Con}(\T_{\Pset}^{\mathrm{loc}})$ would have exactly the additional strength needed to derive $\operatorname{Con}(\T_{\Pset}^{\mathrm{loc}})$ in the local theory.
\end{remark}

\begin{proposition}[Relative consistency step]
\label{prop:relative-consistency-step}
Let $\T_i$ be a recursively enumerable extension of $\mathsf{HA}$, let $C_i$ be a closed sentence of its recursively presented language, and suppose
\begin{equation}
	\T_i\vdash\bigl(C_i\to\operatorname{Con}(\T_i)\bigr).
\end{equation}
Put
\begin{equation}
	\T_{i+1}=\T_i+C_i.
\end{equation}
Then
\begin{equation}
	\T_{i+1}\vdash\operatorname{Con}(\T_i).
\end{equation}
If $\T_{i+1}$ is consistent and carries its standard provability predicate, then
\begin{equation}
	\T_{i+1}\nvdash\operatorname{Con}(\T_{i+1}).
\end{equation}
\end{proposition}

\begin{proof}
	The extension proves $C_i$ as an axiom and inherits the reduction from $\T_i$, so it proves $\operatorname{Con}(\T_i)$. Adding one sentence preserves recursive enumerability and yields another extension of $\mathsf{HA}$. Under the additional consistency hypothesis, Theorem~\ref{thm:conditional-internal-reflection} applied to $\T_{i+1}$ gives the second conclusion.
\end{proof}

\begin{remark}
	Consistency of $\T_i+C_i$ requires a separate metatheoretic assumption at every successor stage, since $\T_i$ may prove $\neg C_i$ even when $\T_i$ is consistent and leaves $C_i$ unprovable. At each accepted stage,
	\begin{equation}
		\T_{i+1}\vdash\operatorname{Con}(\T_i),\quad \T_{i+1}\nvdash\operatorname{Con}(\T_{i+1}).
	\end{equation}
	The proof predicate changes with the theory, so $\operatorname{Con}(\T_i)$ and $\operatorname{Con}(\T_{i+1})$ are different arithmetical assertions. This finite successor pattern is a restricted instance of the proof-theoretic progression framework studied by \citet{feferman62}.
\end{remark}

\begin{proposition}[Conditional L\"obian closure block]
	\label{prop:conditional-loebian-block}
	Let $\T$ satisfy the hypotheses of Theorem~\ref{thm:conditional-internal-reflection}, and let $C$ be a closed sentence of the recursively presented language of $\T$ such that
	\begin{equation}
		\T\vdash\bigl(C\to\operatorname{Con}(\T)\bigr).
	\end{equation}
	Hence,
	\begin{equation}
		\T \nvdash \bigl(
		\Prov_{\T}(\ulcorner C\urcorner)\to C
		\bigr).
	\end{equation}
\end{proposition}

\begin{proof}
	If $\T\vdash\bigl(\Prov_{\T}(\ulcorner C\urcorner)\to C\bigr)$, the arithmetical \textsc{L\"ob Theorem} gives $\T\vdash C$ \citep{loeb55}. This contradicts Proposition~\ref{prop:consistency-reduction-obstruction}.
\end{proof}

\begin{remark}
	Proposition~\ref{prop:unit-removal-trilemma} shows that a prime-related sentence first presupposes a regulated divisibility and primality regime. Prime-power coding then supplies a natural arithmetization of finite sequences, formulas, and derivations. An explicit internal reduction gives such a sentence $C$ the role of a relative consistency premise:
	\begin{equation}
	\T\vdash\bigl(C\to\operatorname{Con}(\T)\bigr).
	\end{equation}
	The arithmetic content of \textsc{Euclidean Escape}, twin-prime recurrence, and \textsc{Goldbach} coverage leaves this reduction as separate input. Propositions~\ref{prop:consistency-reduction-obstruction} and~\ref{prop:conditional-loebian-block} use an internal $\T$-derivation of the reduction; an external derivation establishes a different relative result.

	The derivation $\T+C\vdash\operatorname{Con}(\T)$ is a rule-conformity fact. Its use as a reliable consistency certificate is relative to a metatheory that accepts the interpretation of $C$, the proof predicate, and the reduction. The progression $\T_i\mapsto \T_i+C_i$ tracks changes in proof-theoretic strength, while the verifier sequence in Remark~\ref{rem:verifier-sequence} tracks semantic dependence. Keeping these two forms of relativity distinct prevents prime closure from being identified with consistency by analogy alone.
\end{remark}

\begin{thesis}[Proof regulation and prime agreement]
\label{thesis:proof-regulation-agreement}
	The ternary classifier stabilizes three admissible forms of arithmetic evidence and supports mechanical checking of proof objects formed under them. Relative to the standard natural numbers, the regulated predicates have their familiar extensions. A proof of agreement with primes relates the formal package to an interpreted arithmetic structure and therefore has the relative form
\begin{equation}
	S\vdash\StdDefAgr_I(\Theta_{\mathrm{U,P,C}}).
\end{equation}
	Internal multiplicative agreement is automatic in a definitional expansion of a fixed model, while the representation of $\N$ and its interpretation as intended arithmetic occur in the metatheory. An explicit internal reduction turns a prime closure principle into a relative consistency premise, and adjoining that premise changes the theory whose consistency is next at issue. \textsc{Rice's Theorem} and \textsc{G\"odel's Second Incompleteness Theorem} impose further, separately scoped limits on uniform machine classification and internal reflection.
\end{thesis}

\section{Certification Regimes}
\label{sec:compression-oraclehood}

\subsection{Recursive enumerability}
\label{sec:compression-oraclehood:witnesses}

\begin{definition}[Witness presentation]
	Let $A\subseteq\N$ have a recursively enumerable presentation
	\begin{equation}
		A(n)\iff\exists w\,Q(n,w),
	\end{equation}
	where $Q$ is decidable. A witness $w$ is a finite positive certificate: enumerating pairs $(n,w)$ and checking $Q(n,w)$ enumerates $A$. Compositehood has this form,
	\begin{equation}
		\Comp(n)\iff
		\exists a,b\,\bigl(a>1\wedge b>1\wedge ab=n\bigr).
	\end{equation}
\end{definition}

\begin{remark}
	For compositehood, the witness presentation is bounded by the input, as in Definition~\ref{def:arithmetic-predicates}, and therefore gives a decidable example of recursive enumerability.
\end{remark}

\subsection{Bilateral-completeness obstruction}
\label{sec:compression-oraclehood:bilateral}

\begin{definition}[Instance set]
	Let $\T$ be effectively presented with decidable proof checking, and let $\alpha(x)$ be interpreted over the standard natural numbers. Write
	\begin{equation}
		A_\alpha=\{n\in\N:\N\models\alpha(\bar n)\}.
	\end{equation}
\end{definition}

\begin{proposition}[Bilateral completeness obstruction]
	\label{prop:bilateral-completeness-obstruction}
	Suppose that $\T$ is sound for all instance sentences $\alpha(\bar n)$ and $\neg\alpha(\bar n)$. If, for every $n$, $\T$ proves one of these two sentences, then $A_\alpha$ is recursive.
\end{proposition}

\begin{proof}
	On input $n$, enumerate $\T$-proofs and inspect their conclusions. The bilateral hypothesis ensures that a proof of one of the displayed instances eventually occurs. Soundness identifies its conclusion with the truth value of $\alpha(\bar n)$, so the procedure decides $A_\alpha$.
\end{proof}

\begin{corollary}[No effective bilateral certification]
	\label{cor:no-effective-bilateral-certification}
	If $A_\alpha$ is recursively enumerable and nonrecursive, no sound effectively presented theory proves, for every $n$, one of $\alpha(\bar n)$ and $\neg\alpha(\bar n)$.
\end{corollary}

\begin{remark}
	Positive certificates can still enumerate $A_\alpha$. The obstruction concerns complete effective production of certificates for both sides: parallel enumeration of sound positive and negative certificates would decide the set.
\end{remark}

\begin{definition}[Bilateral proof length]
	For a fixed proof-size measure, define
	\begin{equation}
		L_{\T}^{\pm}(n)=\min\left\{\lvert d\rvert:
		\operatorname{Prf}_{\T}(d,\ulcorner\alpha(\bar n)\urcorner)
		\vee
		\operatorname{Prf}_{\T}(d,\ulcorner\neg\alpha(\bar n)\urcorner)\right\},
	\end{equation}
	with value $\infty$ when neither proof exists.
\end{definition}

\begin{proposition}[Recursive bilateral bounds decide]
	\label{prop:recursive-bilateral-bound}
	Assume that $\T$ is sound for the displayed instance formulas. If a total recursive function $g$ satisfies $L_{\T}^{\pm}(n)\leq g(n)$ for every $n$, then $A_\alpha$ is recursive.
\end{proposition}

\begin{proof}
	For input $n$, check the finitely many codes of size at most $g(n)$. One is a proof of one of the two instance formulas; soundness identifies the correct side.
\end{proof}

\begin{corollary}[Tagged certificate producer computes the classified set]
	\label{cor:tagged-compressor-oracle}
	Let $X$ be an oracle. Suppose an $X$-recursive total function returns $F^X(n)=(b_n,d_n)$, where $b_n\in\{0,1\}$ and
	\begin{equation}
		b_n=1\Rightarrow\operatorname{Prf}_{\T}(d_n,\ulcorner\alpha(\bar n)\urcorner),\quad b_n=0\Rightarrow\operatorname{Prf}_{\T}(d_n,\ulcorner\neg\alpha(\bar n)\urcorner).
	\end{equation}
	If $\T$ is sound for these formulas, then $A_\alpha\leq_T X$.
\end{corollary}

\begin{proof}
	The tag $b_n$ decides whether $n\in A_\alpha$.
\end{proof}

\begin{remark}
	Under \textsc{Curry--Howard}, the analogous type is
	\begin{equation}
		\prod_{n:\N}\bigl(\alpha(n)+\neg\alpha(n)\bigr).
	\end{equation}
	In a type theory with effective normalization and canonicity for sums, an inhabitant computes a branch for each numeral. An opaque axiom of this type supplies abstract inhabitation without an executable classifier; operational branch exposure is the additional requirement at issue here \citep{howard80,martinlof84,sorensen06}.
\end{remark}

\subsection{\textsc{Busy Beaver} as a completion bound}
\label{sec:compression-oraclehood:busy-beaver}

\begin{definition}[Running-time \textsc{Busy Beaver} function]
	Fix a deterministic acceptable universal programming system $U$: its step relation is primitive recursive and it effectively simulates every partial recursive program. In particular, $K_U$ is a complete recursively enumerable set. Let $\tau_U(p)$ be the running time of a halting program $p$, and define
	\begin{equation}
		\operatorname{BB}_U(n)=
		\max\bigl(\{\tau_U(p):\lvert p\rvert\leq n\wedge U(p)\downarrow\}\cup\{0\}\bigr).
	\end{equation}
	Put $K_U=\{p:U(p)\downarrow\}$. For such a system, the finite maximum is nonrecursive \citep{rado62,rogers67}.
\end{definition}

\begin{theorem}[\textsc{Busy Beaver}--oracle equivalence]
	\label{thm:busy-beaver-oracle-equivalence}
	The halting set computes $\operatorname{BB}_U$. Conversely, any oracle $X$ computing a total majorant $b$ with $b(n)\geq\operatorname{BB}_U(n)$ for all $n$ computes $K_U$. Hence
	\begin{equation}
		\operatorname{BB}_U\equiv_T K_U\equiv_T 0'.
	\end{equation}
\end{theorem}

\begin{proof}
	Given $K_U$ and $n$, identify the finitely many halting programs of length at most $n$, run them to termination, and take the largest running time. Conversely, on input $p$, compute $b(\lvert p\rvert)$ and simulate $U(p)$ for that many steps. Failure to halt by then proves that it never halts.
\end{proof}

\begin{definition}[Completion formula]
	Let $\operatorname{FirstHalt}_U(p,t)$ be the primitive-recursive relation saying that $U(p)$ first halts at step $t$, and put
	\begin{equation}
		\Theta_U(n,B)\equiv
		\forall p,t\,\bigl(\lvert p\rvert\leq n\wedge\operatorname{FirstHalt}_U(p,t)\to t\leq B\bigr).
	\end{equation}
	Over the standard natural numbers, $\Theta_U(n,B)$ holds exactly when $B\geq\operatorname{BB}_U(n)$.
\end{definition}

\begin{proposition}[Constructive completion produces an oracle]
	\label{prop:constructive-completion-oracle}
	Let $\T$ be effectively presented and sound for the formulas $\Theta_U(\bar n,\bar B)$. If an oracle $X$ computes pairs $(B_n,d_n)$ with
	\begin{equation}
		\operatorname{Prf}_{\T}\bigl(d_n,\ulcorner\Theta_U(\bar n,\overline{B_n})\urcorner\bigr)
	\end{equation}
	for every $n$, then $0'\leq_T X$.
\end{proposition}

\begin{proof}
	Soundness makes $n\mapsto B_n$ an $X$-recursive \textsc{Busy Beaver} majorant. Theorem~\ref{thm:busy-beaver-oracle-equivalence} then gives $K_U\leq_T X$.
\end{proof}

\begin{corollary}
	No total recursive function uniformly produces sound completion certificates of the displayed form.
\end{corollary}

\begin{proof}
	Such a function would instantiate Proposition~\ref{prop:constructive-completion-oracle} with the empty oracle and compute $0'$ recursively.
\end{proof}

\begin{definition}[Trace and oracle calculi]
	Let $P_{\mathrm{tr}}$ be a trace calculus in which a proof of $U(\bar p)\downarrow$ explicitly lists every configuration from the initial state to the halting state. Assume the chosen encoding gives a constant $c>0$ such that a computation of $t$ steps requires a proof of size at least $ct$. Extend it to $P_{\mathrm{tr}}^K$ with the oracle rules
	\begin{equation}
		\frac{}{U(\bar p)\downarrow}\;(p\in K_U),\quad\frac{}{\neg U(\bar p)\downarrow}\;(p\notin K_U).
	\end{equation}
	Rule conformity is decidable relative to $K_U$ \citep{turing37}.
\end{definition}

\begin{theorem}[BB speed-up]
	\label{thm:oracle-busy-beaver-speedup}
	With binary program numerals, there are constants $c,C>0$ and, for all sufficiently large $n$, halting programs $p_n$ of length at most $n$ such that
	\begin{equation}
		L_{P_{\mathrm{tr}}}\bigl(U(\overline{p_n})\downarrow\bigr)\geq c\operatorname{BB}_U(n),\quad L_{P_{\mathrm{tr}}^K}\bigl(U(\overline{p_n})\downarrow\bigr)\leq Cn.
	\end{equation}
	No total recursive function bounds all $P_{\mathrm{tr}}$ proof lengths by a function of the corresponding $P_{\mathrm{tr}}^K$ proof lengths.
\end{theorem}

\begin{proof}
	Choose a program attaining $\operatorname{BB}_U(n)$. Its trace has the stated lower bound; one oracle-rule instance with a binary numeral has size $O(n)$. If a recursive comparison $f$ existed, then $\max_{m\leq Cn}f(m)$ would be a recursive upper bound for the corresponding trace lengths and hence for $c\operatorname{BB}_U(n)$. This contradicts Theorem~\ref{thm:busy-beaver-oracle-equivalence}.
\end{proof}

\begin{remark}
	The converse needs bilateral, branch-exposing, or completion information. Short positive halting traces alone leave nonhalting undecided.
\end{remark}

\subsection{Relativization}
\label{sec:compression-oraclehood:relativization}

\begin{proposition}[Relativized completion boundary]
	For every oracle $X$ and universal $X$-oracle \textsc{Turing Machine} $U^X$,
	\begin{equation}
		X\oplus\operatorname{BB}_{U}^{X}\equiv_T X',
	\end{equation}
	where $X'$ is the $X$-relative halting problem. A pointwise $X$-relative completion bound computes the next jump.
\end{proposition}

\begin{proof}
	The $X$-jump identifies the halting $X$-programs of a fixed bounded length, after which $X$ simulates them to obtain the maximum running time. Conversely, $X$ together with any such majorant decides $X$-program halting by bounded simulation.
\end{proof}

\begin{remark}
	The local predicates $\Comp(n)$ and $\Prime(n)$ remain at the decidable base: their factor spaces are bounded by $n$. If $G(n)$ says that $2n$ has a \textsc{Goldbach} representation, then $G$ is decidable by bounded search and \textsc{Goldbach} coverage is a $\Pi^0_1$ sentence. Twin-prime recurrence,
	\begin{equation}
		\forall B\,\exists p>B\,\bigl(\Prime(p)\wedge\Prime(p+2)\bigr),
	\end{equation}
	has $\Pi^0_2$ form. These descriptions concern uniform indexed families: a fixed conjecture has one truth value, while oracle content belongs to a uniform operation deciding instances or returning completion bounds.

	Adding one conjecture as an axiom preserves recursive enumerability. Oracle strength arises from uniform operational content, such as a map from each instance to a correct tagged proof or to a sound completion bound.
\end{remark}

\begin{thesis}[Uniform completion boundary]
	\label{thesis:compression-oracle-boundary}
	Finite witness presentations supply one-sided evidence. Proof-regime speed-up may shorten derivations without deciding a nonrecursive family. A uniform, sound operation that exposes the correct branch for every instance, or that closes universal halting search by finite completion bounds, computes the classified set. For universal computation, completion bounds are equivalent in degree to a \textsc{Busy Beaver} majorant and hence to the halting oracle; relativized completion yields the corresponding oracle jump. None of these obstructions applies to the bounded prime classifier on $\N$.
\end{thesis}

\begin{thesis}[Primes as mirage]
	\label{thesis:prime-holes-global-regularity}
	Relative to interior multiplication, the composites form the positively generated image
	\begin{equation}
		\operatorname{im}(\mu_I),
	\end{equation}
	while the primes are the holes in that image:
	\begin{equation}
		\NgtOne\setminus\operatorname{im}(\mu_I).
	\end{equation}
	A composite is certified by a preimage under multiplication; a prime is certified by the absence of every such. This asymmetry persists even though each bounded instance is decidable.

	Global prime conjectures ask for recurrence, pairing, or additive organization among these holes without supplying a corresponding positive generator for their distribution. Arbitrarily large finite regions may therefore be verified while the unbounded pattern remains unproved. In this proof-polarity sense, the co-recursively enumerable presentation of the primes records absence of construction.
\end{thesis}

\clearpage

\section{Conclusion}
\label{sec:conclusion}

Over $\mathbb{N}$, unit, prime, and composite form a constructive ternary classification whose evidence consists respectively of equality, an interior factorization, and a bounded family of refutations. Because the candidate factor domain for each input is finite and decidable, this classifier is uniformly executable. The staged formulation exposes the data on which that execution depends: a presented block family, an admissible order, and fresh compatibility information under extension. Completion by the bounded factor square converts stage-relative survival into primality.

The arithmetic developments show how this classification propagates. Every finite catcher is sound while missing some composites; adjoining each least survivor as a new generator recursively enumerates the primes, and the square of the next prime is the least composite missed at that stage. Algebraically, the same prime--composite boundary appears as the contrast between modular cancellation and explicit nonzero annihilation. Prime factorization, coordinate normal forms, and finite cancellation certificates extend this boundary into effective proof objects.

\begin{figure}[ht]
	\centering
	\includegraphics[width=.65\textwidth]{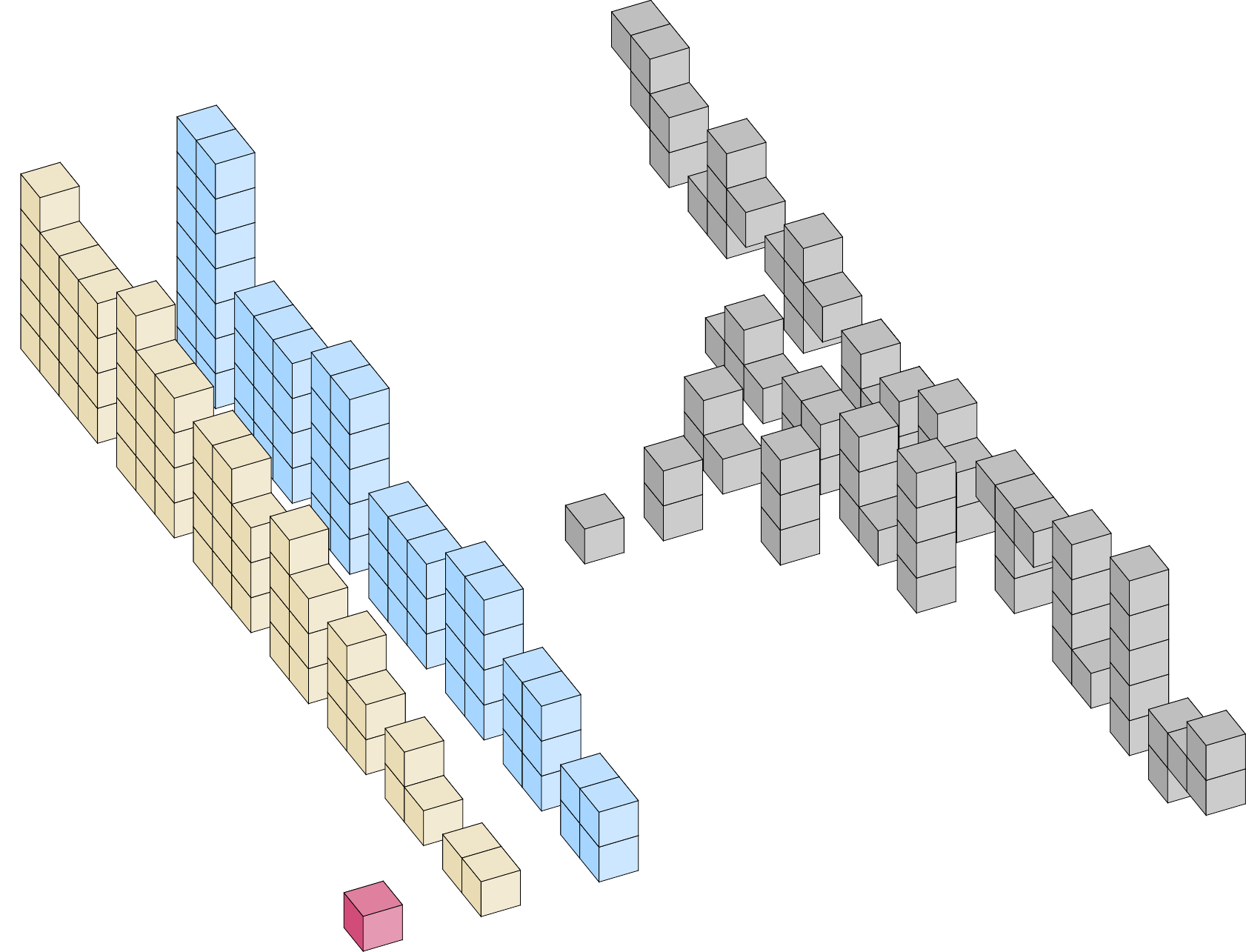}
	\caption{Internal satisfaction of the defining clauses fixes the predicate extensions relative to a model. Representation of $\mathbb{N}$ and the interpretation of definitional agreement as agreement with intended arithmetic occur in a (colored) metatheory; consistency and reflection occupy further, distinct levels.}
	\label{fig:semantic-metatheory}
\end{figure}

The semantic and computational results locate these constructions within their supporting frameworks. The defining clauses yield a conservative extension of $\mathsf{HA}$ and a unique expansion of every arithmetic model. Their agreement is therefore relative to the represented structure, as summarized in Figure~\ref{fig:semantic-metatheory}; standardness, soundness, consistency, and reflection require further metatheoretic data. Uniform bilateral certification and universal completion bounds carry additional computational strength, reaching the halting degree and its relativized jumps. The bounded prime classifier remains decidable because each input supplies its own finite search space.

Taken together, these results separate three levels that constructive arithmetic must coordinate: local evidence for a numerical classification, uniform organization of that evidence, and semantic interpretation of the resulting formal theory.
		\clearpage
\section*{References and Notes}

	{\scriptsize
		\bibliographystyle{plainnat}
		\setlength{\bibsep}{0.1em}
		\bibliography{refs}}

		\clearpage
\vspace*{\fill}

\subsection*{Note on the Third Version}

This third version corrects the overextended framing of the second. The earlier manuscript treated bounded factor search as an unbounded hierarchy and connected primality too directly with oracle and global-conjecture claims. Here factor search over the standard natural numbers remains bounded and decidable. The paper studies the proof-theoretic regulation of unit, composite, and prime, and the relative conditions under which formal rules receive an arithmetic interpretation.

\subsection*{Final Remarks}

The author welcomes criticism, proposed extensions, scholarly correspondence, and constructive dialogue. No conflicts of interest are declared. This research received no funding.

\begin{center}
	\scriptsize{
	\vspace{1em}
	Milan Rosko\\
	\vspace{1em}
	ORCID: \href{https://orcid.org/0009-0003-1363-7158}{\textsf{0009-0003-1363-7158}}\\[1ex]
	Email: \href{mailto:hi-at-milanrosko.com}{\textsf{hi-at-milanrosko.com}}\\[1ex]
	Email: \href{mailto:Q1012878@studium.fernuni-hagen.de}{\textsf{Q1012878@studium.fernuni-hagen.de}}\\
	\vspace{1em}
	Licensed under {\ccby}\\ \href{https://creativecommons.org/licenses/by/4.0/}{\scriptsize\textsf{creativecommons.org/licenses/by/4.0}}
	}
\end{center}
\vspace*{\fill}